\documentclass[10pt,journal,compsoc]{IEEEtran}
\usepackage{ctable}

\ifCLASSOPTIONcompsoc
  \usepackage[nocompress]{cite}
\else
  \usepackage{cite}
\fi
\ifCLASSINFOpdf
\else
\fi
\usepackage{pdfpages}
\usepackage{graphicx}
\usepackage{float}
\usepackage{amssymb,amsfonts,amsmath,amsthm}
\newcommand{\R}{\mathbb{R}} 
\newtheorem{theorem}{Theorem}
\newtheorem{lemma}{Lemma}

\usepackage{breqn}
\usepackage{algorithm,algpseudocode}
\usepackage{enumitem}
\begin{document}
%
\title{Constrained Mass Optimal Transport}
%
%
%
%

\author{Said~Kerrache,~\IEEEmembership{Member,~IEEE,}
        and~Yasushi~Nakauchi,~\IEEEmembership{Member,~IEEE}
\IEEEcompsocitemizethanks{\IEEEcompsocthanksitem S. Kerrache was with the Graduate School of Systems and Information Engineering, University of Tsukuba, 1-1-1 Tennodai, Tsukuba, Ibaraki 305-8577 Japan.\protect\\
E-mail: see http://www.michaelshell.org/contact.html
\IEEEcompsocthanksitem Y. Nakauchi is with the Graduate School of Systems and Information Engineering, University of Tsukuba, 1-1-1 Tennodai, Tsukuba, Ibaraki 305-8577 Japan}
\thanks{Manuscript received April 19, 2005; revised August 26, 2015.}}

%
%

\markboth{Journal of \LaTeX\ Class Files,~Vol.~14, No.~8, August~2015}%
{Shell \MakeLowercase{\textit{et al.}}: Bare Demo of IEEEtran.cls for Computer Society Journals}
%



\IEEEtitleabstractindextext{%
\begin{abstract}
Optimal mass transport, also known as the earth mover's problem, is an optimization problem with important applications in various disciplines, including economics, probability theory, fluid dynamics, cosmology and geophysics to cite a few.  Optimal transport has also found successful applications in image registration, content-based image retrieval, and more generally in pattern recognition and machine learning as a way to measure dissimilarity among data. This paper introduces the problem of constrained optimal transport. 
The time-dependent formulation, more precisely, the fluid dynamics approach is used as a starting point from which the constrained problem is defined by imposing a soft constraint on the density and momentum fields or restricting them to a subset of curves that satisfy some prescribed conditions. 
A family of algorithms is introduced to solve a class of constrained saddle point problems, which has convexly constrained optimal transport on closed convex subsets of the Euclidean space as a special case. Convergence proofs and numerical results are presented. \end{abstract}

}

\maketitle

\IEEEdisplaynontitleabstractindextext

%
\IEEEpeerreviewmaketitle

\IEEEraisesectionheading{\section{Introduction}\label{sec:introduction}}

Optimal transport is a mathematical problem with a long and rich history.
It has classically been found to have important applications in economics \cite{kantorovitch42} and probability theory \cite{rachev98}.
In recent years, the research on the subject has undergone a rapid expansion at  the theoretical as well as the application level. 
Optimal transport has been applied and connected to various fields in mathematics, science and engineering, including geometry \cite{mccan10}, fluid dynamics \cite{jordan98,carrillo07}, optics \cite{gutierrez09,rubinstein07}, geophysics and oceanography \cite{cullen06}, ecology \cite{Kranstauber17}, meteorology \cite{cullen03}, cosmology \cite{frisch02,brenier03}, antenna design \cite{glim03}, image processing and computer vision \cite{gangbo00,haker04,museyko09}, information retrieval and pattern recognition \cite{jang18,Hou2016SquaredEM,xu16,Fan2018,Yuan2018,Atasu2019, Zhang16, Charles17, Calarasanu15,Wang19,Uysal15}. Our aim in this work is to solve the optimal mass transportation problem under constraints that could be either hard, soft or a combination of both.

Applications of optimal transport can be classified into two major categories.
Applications in the first category have as primary target finding the optimal transport plan or the displacement interpolation.
This is the case for problems of physical nature, such as planning optimal transport routes for goods or tracking back the universe in time \cite{frisch02}.
The possibility of solving constrained transport problems undoubtedly opens new areas of applications within this category.
It allows, for instance, to focus on a subset of mass distributions of interest, or to eliminate those distributions that violate some given laws, physical or otherwise.

Applications in the second category seek to compute the distance between densities.
This is typically the case for recognition and matching problems.
The densities considered in this type of applications are not arbitrary.
They represent instances of a certain family of patterns, for example, images of human faces or medical images of a certain body organ.
Constraining the interpolation to the set of densities representing the family of patterns under consideration allows to compute the distance induced by the ambient optimal transport metric on the feasible set, thus probing its intrinsic geometry.
The importance of intrinsic geometry for revealing the true structure of  data has come to be realized by the researchers in the field of pattern recognition and data analysis and was at the root of the emergence of the field of nonlinear dimensionality reduction \cite{isomap,lle}. This work introduces numerical algorithms that can be used to practically compute the optimal mass transport interpolation under constraints and the associated distance.

The content of this paper is organized as follows.
Related work is reviewed in Section \ref{sec:related-work}. Section \ref{sec:problem} introduces the problem formally using the fluid dynamic formulation.
In Section \ref{sec:computational-methods}, a family of numerical algorithms for solving convexly constrained transport problems with the squared Euclidean distance as cost is presented.
Numerical experiments demonstrating the working and effectiveness of the algorithms ar presented in Section \ref{sec:experiments}.
The final section concludes the paper with a discussion of a number of aspects related to this work.

\section{Related Work}
\label{sec:related-work}

The time-dependent formulation of optimal transport, otherwise known as displacement interpolation, was introduced by McCann \cite{mccann97} in the case of the squared Euclidean distance as a way to model certain gas dynamics.
Benamou and Brenier established that McCann's interpolation can be seen as an action minimization problem in the space of probability measures \cite{benamou00cfd}, an equivalence that turned out to be one of the most important results in the theory of optimal transport.
Otto \cite{otto01} went further in this direction and formulated the equation of Benamou and Brenier as that of a geodesic in a Riemannian setting, where the space of probability measures is interpreted as a manifold with the Wasserestein metric being the associated Riemannian metric.
Villani and Lott \cite{villani09,lott09} generalized displacement interpolation to length spaces, where no smooth structure is required.
This is achieved through the notion of dynamical transference plan and abstract Lagrangian action.

The idea of restricting the set of admissible interpolations was introduced in \cite{brenier-const97,brenier-const99,brenier-const02,brenier-const03} in relation with optimal multiphase transportation.
The authors, however, do not attempt to numerically solve the problem.
It is also customary to restrict the set of admissible interpolations when investigating theoretical questions related to the existence and regularity of the optimal transport solutions.
However, this type of research is not concerned with finding numerical algorithms to solve the resulting problems.
The work of \cite{brenier-const97,brenier-const99,brenier-const02,brenier-const03} inspired \cite{buttazzo} to solve the optimal transport problem penalized with soft constraints.

The earliest practical method for solving the optimal transport problem was the linear programming method proposed by Kantorovich \cite{kantorovitch42}.
This method is based on the relaxation of the optimal transport problem proposed by Kantorovich himself.
The principle of the relaxation is to enlarge the set of possible transport plans to allow for mass splitting, that is, eliminating the requirement of the determinism of the transport map.
This gives rise to a duality result and the formulation of the problem as a linear program.

Similar in nature to Kantorovich's method, a discrete method to compute the optimal transport map is used in \cite{frisch02,brenier03,cosmology08} for the reconstruction of the early state of the universe.
For instance, in \cite{frisch02}, the optimal transport is modeled as an assignment problem and solved using efficient linear programming algorithms.
Note that the underlying assumption in these applications is that matter is discrete in the universe, hence the discrete formulation of the optimal transport problem.
Despite the success of the linear programming model, this method has its shortcomings.
Mainly, it computes only the time-independent solution to the problem, which makes it virtually impossible to impose any constraint on the interpolation.

Within the image and computer vision community, solving the optimal transport problem is almost exclusively accomplished by calculating the optimal transport map \cite{haker01,haker03,haker03-tracking,haker04,haker05,zhu07,museyko09,gpu,dominitz10}.
This practice has its origins in the field of image registration, which remains the predominant area of application of optimal transport in imagery.
More precisely, the goal is to find a map that minimizes the cost
\begin{dmath}
\int \left | u\left ( x \right ) -x \right |^2 \mu_0\left ( x \right )dx,
\end{dmath}
among all mass preserving maps $u$ that transport the initial density $\mu_0$ to the final density $\mu_1$.
Here, $u$ associates to each point $x$ its destination $u\left(x\right)$.
This approach has the merit of directly computing the optimal map, which is very important in the case of image registration.
The optimization methods used in this approach are of gradient descent type, which renders the method more computationally efficient than the computational fluid dynamics approach \cite{gpu}.
On the other hand, it has the disadvantage of not offering a direct control on the intermediate densities, which makes the introduction of constraints a much difficult task.
It is worth mentioning that there exists a strong link between the computation of the optimal transport map and the solution to the Monge-Ampere equation \cite{loeper05,mohammadi07,cosmology08,Sulman11}. 

In \cite{benamou00cfd,benamou01,benamou02,guittet01,guittet03}, the problem of optimal mass transportation in $\R^d$ with the squared Euclidean distance as cost is recast as an optimal control problem of a potential flow.
The approach consists in computing a flow, in the sense of fluid dynamics, that moves the initial density to the final density.
The optimal flow turns out to be potential, thus it is fully controlled by a potential function to be determined by non-linear optimization.
This method not only finds the transport map, but also provides a time interpolant of the density.
The three main steps in this approach are: 
\begin{itemize}
	\item Establishing the equivalence between optimal transport and the minimization of the kinetic energy of the flow.
	\item Converting the non-convex kinetic energy minimization problem into a convex one, through an appropriate change of variables.
	\item Changing the minimization problem into a saddle point problem through the use of the Fenchel transform (or the Legendre transform).
\end{itemize}
The saddle point problem is finally solved using an Uzawa type algorithm described in \cite{fortin83,glowinski89}.
This method is an efficient and elegant way to solve the optimal transport problem.
It provides a direct access to the interpolating densities as well as the velocity field.
This makes it a perfect starting point from which to tackle the constrained problem.
However, the form of the saddle point problem makes the application of constraints a non-trivial task, and one of the aims of this work is to find an efficient solution to this problem.

In \cite{buttazzo}, the authors present a framework that mimics \cite{benamou00cfd} to compute a transport flow that minimizes different types of cost functions, and not only the kinetic energy.
The cost functions presented are related to the problem of transportation in congested domains, a problem with important applications to crowd movement modeling \cite{picolli09,picolli10}.
The framework allows, theoretically, for a wide class of cost functions.
It is possible to penalize certain transport plans by adding a penalty term to the initial cost function, which can be the kinetic energy for instance.

In addition to this soft constraint handling approach, it is theoretically possible to apply hard constraints by adding indicator functions.
In practice, however, the method requires the computation of the Fenchel transform of the resulting cost functions.
This can be done analytically for very simple penalty and characteristic functions.
The authors carried out the transformation for three examples, two of them contain penalty terms, whereas the only example containing an indicator function (a hard constraint) is the case of a simple bound constraint on the density.
Such computation is, however, impossible for complex cases, and resorting to numerical optimization defeats the purpose and the usefulness of the method.
In contrast, the algorithms proposed in this work are not limited to simple constraints and do not require any prior analytic computation.

\section{Problem Statement}
\label{sec:problem}

The computational fluid dynamic  formulation of optimal transport has its roots in the seminal work of Benamou and Brenier \cite{benamou00cfd}.
The authors established the equivalence on $\R^d$ between the kinetic energy functional \eqref{eq:bb-lagrangian-velocity} and the Wasserstein metric $W_2$.
This result lead to the development of a numerical scheme for solving the optimal transport problem with the squared Euclidean distance as cost under appropriate regularity conditions on the density.
Given a closed convex subset $D$ of $\R^d$, on which two densities $\rho_0$ and $\rho_1$ are prescribed, the problem addressed in  \cite{benamou00cfd} is finding a density $\rho$ and a momentum $m$ that solve
\begin{dmath}
\label{eq:bb-lagrangian-velocity}
\inf_{\rho, m} \int_0^1 \int_D\frac{\left | m\left ( t,x \right ) \right |^2}{2\rho\left ( t,x \right )}dx dt,
\end{dmath}
subject to the following continuity equation and boundary conditions
\begin{equation}
\label{eq:continuity-boundary-conditions}
\partial_t \rho + \nabla \cdot m =0,  \quad \rho\left ( 0, \cdot \right )= \rho_0, \quad \rho\left ( 1,\cdot \right )= \rho_1.
\end{equation}
The importance of this numerical scheme to the problem of constrained optimal transport is that it computes not only the transport map but the whole displacement interpolation.
Therefore, solving a constrained problem amounts to adding the penalty term $I \left( \rho, m \right) $.
That is
\begin{dmath}
\label{eq:constrained-lag}
\inf_{\rho,m} \left[\int_0^1\int_D \frac{\left | m \left ( t,x \right ) \right |^2}{2 \rho \left ( t,x \right )} dx dt + I \left (\rho,m \right)\right],
\end{dmath}
subject to \eqref{eq:continuity-boundary-conditions}.
In the hard constraint case, $I \left( \rho, m \right) $ corresponds to the indicator function of a set $U$. Thanks to this formulation, not only can we impose impose constraints on the intermediate density distributions, but also on the momentum field of the flow, hence allowing a finer control on the transportation process. In the rest of this paper, we will introduce numerical algorithms for solving this problem for the case where $ I $ is a convex functional.

\section{Computational Methods}
\label{sec:computational-methods}
Following \cite{benamou00cfd} and using the same notation, Problem \eqref{eq:constrained-lag} can be written as 
\begin{dmath}
\label{eq:original-lagrangian}
\inf_{\phi,q} \sup_{\mu} F\left ( q \right )+G\left ( \phi \right ) - I \left (\mu \right) +\left \langle \mu, \nabla \phi -q \right \rangle,
\end{dmath}
where $\mu= \left ( \rho, m \right)$, $\phi$ is the Lagrange multiplier of the constraint \eqref{eq:continuity-boundary-conditions}, $G\left ( \phi \right )=\int_D \phi\left ( 0,x \right )\rho_0\left ( x \right )-\phi\left ( 1,x \right )\rho_1\left ( x \right ) dx$, $F$ is the indicator function for the set  
\begin{multline}
K= \{ \left ( a,b \right ): \R \times \R^d \to \R \times \R^d, \\a+\left | b \right |^2/2 \leq 0 \text{ pointwise}  \},
\end{multline}
and $\left \langle \cdot, \cdot \right \rangle$ is the inner product defined by $\left \langle u,v \right \rangle=\int_0^1 \int_D u \cdot v $.

The algorithms presented in this section solve a general class of problems, of which Problem \eqref{eq:original-lagrangian} is a special case.
Given two real Hilbert spaces $V$ and $H$, with inner products denoted by $\left \langle \cdot , \cdot \right \rangle$ and norms denoted by $\left | \cdot  \right |$, the class of problems considered is 
\begin{multline}
\label{eq:general-lagrangian}
\inf_{\phi \in V,q \in H} \sup_{\mu \in H}\mathcal{L}\left ( \phi,q,\mu \right )= F\left ( q \right )+G\left ( \phi \right ) - I \left ( \mu \right ) \\ + \left \langle \mu, B \phi -q \right \rangle,
\end{multline}
where $F$, $I$ and $G$ are convex, lower semi-continuous functionals on $H$ and $V$ respectively and $B: V \to H$ is a linear operator.
In the case where $I$ is an indicator function, the convexity and the lower semi-continuity translates into the feasible set $U$ being convex and closed respectively.
The additional term $I \left ( \mu \right )$ prevents from using the same augmented Lagrangian technique as in \cite{fortin83,glowinski89}, since it is possible to have
\begin{dmath}
B \phi^* - q^* \neq 0
\end{dmath}
at a saddle point $\left (\phi^*,q^*,\mu^* \right)$ of $\mathcal{L}$.
We will show, however, that the problem can be transformed to a computationally manageable form by an appropriate use of the decomposition-coordination principle.
The first step to this end is to decouple $\mu$ from $\phi$ and $q$.
This is achieved by introducing two additional variables $p$ and $b$, and imposing the constraints $B \phi = p$ and $b=q$.
Problem \eqref{eq:general-lagrangian} then takes the form 
\begin{multline}
\label{eq:lagrangian-p-b}
\inf_{\phi, q, p, b} \sup_{\mu, \nu, \eta}\mathcal{L}'\left ( \phi,q,p,b,\mu, \nu,\eta \right )=  F\left ( q \right )+G\left ( \phi \right ) \\ - I \left ( \mu \right )
+\left \langle \mu, p - b \right \rangle + \left \langle \nu , B \phi - p \right \rangle + \left \langle \eta , b - q \right \rangle, 
\end{multline}
where $\nu$ and $\eta$ are the Lagrange multipliers for the new constraints.
A key observation at this point is that any saddle point $\left ( \phi^*, q^*, p^*, b^*, \mu^*, \nu^*, \eta^* \right )$ of Problem \eqref{eq:lagrangian-p-b} verifies $\mu^*=\nu^*=\eta^*$ as shown in Theorem \ref{thm:lagrangian-properties}.
This relation allows to augment the Lagrangian by  providing additional concavity, which proves crucial for solving the problem.
The new form upon augmentation is 
\begin{dmath}
\label{eq:aug-lagrangian}
\inf_{\phi, q, p, b} \sup_{\mu, \nu, \eta} \mathcal{L}_{r,s}\left ( \phi,q,p,b,\mu, \nu,\eta \right )= F\left ( q \right )+G\left ( \phi \right )-I \left ( \mu \right )
+\left \langle \mu, p - b \right \rangle + \left \langle \nu , B \phi - p \right \rangle + \left \langle \eta , b - q \right \rangle 
+ \frac{r}{2} \left | B \phi - p \right |^2 + \frac{r}{2} \left | b - q \right |^2 - \frac{s}{2} \left | \mu - \nu \right |^2 - \frac{s}{2} \left | \mu - \eta \right |^2.
\end{dmath}
The three problems \eqref{eq:general-lagrangian}, \eqref{eq:lagrangian-p-b} and \eqref{eq:aug-lagrangian} are equivalent in the sense that they essentially possess the same set of saddle points.
This is the content of Theorem \ref{thm:lagrangian-properties}.
\begin{theorem}[Equivalence of problems  \eqref{eq:general-lagrangian}, \eqref{eq:lagrangian-p-b} and \eqref{eq:aug-lagrangian}]
	\label{thm:lagrangian-properties}
	Assume that $G$, $F$, $I$ are proper convex functional on $V$, $H$ respectively and $B$ is a linear map from $V$ to $H$.
	Let $r$, $s>0$ and $\left ( \phi^*, q^*, p^*, b^*, \mu^*, \nu^*, \eta^* \right ) \in V \times H^6$. 
	Then, the following statements are equivalent:
	\begin{enumerate}
		\item $\left ( \phi^*, q^*,  \mu^* \right )$ is a saddle point of $\mathcal{L}$, $p^*= B \phi^*$, $b^*=q^*$ and $\mu^*=\nu^*=\eta^*$. \item $\left ( \phi^*, q^*, p^*, b^*, \mu^*, \nu^*, \eta^* \right) $ is a saddle point of  $\mathcal{L}'$.
		\item $\left ( \phi^*, q^*, p^*, b^*, \mu^*, \nu^*, \eta^* \right) $ is a saddle point of $\mathcal{L}_{r,s}$.
	\end{enumerate}
\end{theorem}
\begin{proof}
	in supporting material.
\end{proof}
Notice that Theorem \ref{thm:lagrangian-properties} assumes the existence of a saddle point for Problem \eqref{eq:general-lagrangian}.
Proving the existence, on the other hand, has to be done on a case-by-case basis, usually using a separation theorem.
Problem \eqref{eq:aug-lagrangian} can be solved using a number of Uzawa type algorithms presented hereafter.
To simplify the algorithmic notation, define for fixed $\mu$, $p$ and $b$ the two Lagrangians $\mathcal{L}_{p, \mu}$ and $\mathcal{L}_{b, \mu}$ by
\begin{multline}
\mathcal{L}_{p, \mu} \left ( \phi, \nu \right )= G\left ( \phi \right ) + \left \langle \nu , B \phi - p \right \rangle \\+ \frac{r}{2} \left | B\phi -p \right |^2 - \frac{s}{2} \left | \nu - \mu \right |^2,
\end{multline}
\begin{multline}
\mathcal{L}_{b, \mu}\left ( q, \eta \right )= F\left( q \right ) + \left \langle \eta , b -q \right \rangle \\+ \frac{r}{2} \left | q -b \right |^2 - \frac{s}{2} \left | \eta - \mu \right |^2.
\end{multline}
Notice that
\begin{dmath}
\mathcal{L}_{r,s}\left ( \phi,q,p,b,\mu, \nu,\eta \right )= \mathcal{L}_{p, \mu} \left ( \phi, \nu \right ) + \mathcal{L}_{b, \mu}\left ( q, \eta \right ) +\left \langle \mu, p - b \right \rangle - I\left( \mu \right).
\end{dmath}
With this notation, the first algorithm for solving Problem \eqref{eq:aug-lagrangian} is Algorithm \ref{alg:1}.

\begin{algorithm*}[!t]
	\begin{algorithmic}
		\State \textbf{Initialization}: $\left (  p^0, b^0,\mu^0 \right ) \in H^3 $ given arbitrarily.
		\Repeat
		\State (1) Compute $ \phi^n, \nu^n $ such that, $ \forall \left ( \phi, \nu \right ) \in V \times H:	\mathcal{L}_{p^n,\mu^n} \left ( \phi^n, \nu \right ) \leq \mathcal{L}_{p^n,\mu^n} \left ( \phi^n, \nu^n \right )  \leq \mathcal{L}_{p^n,\mu^n} \left ( \phi, \nu^n \right ) $;\\
		\State (2) $ p^{n+1} = p^n - \rho \left ( r + 1/s \right )\left ( p^n - B\phi^n \right ) $; \\
		\State (3) Compute $ q^n $, $ \eta^n $  such that, $ \forall \left ( q, \eta \right ) \in H \times H: \mathcal{L}_{b^n,\mu^n} \left ( q^n, \eta \right ) \leq \mathcal{L}_{b^n,\mu^n} \left ( q^n, \eta^n \right ) \leq \mathcal{L}_{b^n,\mu^n} \left ( q, \eta^n \right ) $;\\
		\State (4) $ b^{n+1} = b^n - \rho \left ( r + 1/s \right )\left ( b^n - q^n \right ) $;\\
		\State (5) $ \mu^{n+1}= \arg\inf_{\mu} [ s | \mu  -   ( \nu^n +\eta^n + ( p^{n+1} - b^{n+1}  )/s )/2 |^2 + I\left ( \mu \right ) ] $;
		\Until{Convergence}
	\end{algorithmic}
	\caption{For solving Problem \eqref{eq:aug-lagrangian}}
	\label{alg:1}
\end{algorithm*}

The sub-problems in Step 1 and 3 are saddle point problems that can be solved using Algorithm \ref{alg:1-1} and \ref{alg:1-2} respectively.
Step 2 and 4 are gradient descent update rules for $p$ and $b$.
The algorithm parameters $\rho$, $\rho_\nu$ and $\rho_\eta$  control the step size for the main and sub-iterations.
Finally, Step 5 is a minimization problem, which reduces to a projection on the feasible set when $I$ is an indicator function.
In the unconstrained case, it reduces further to the simple update:
\begin{dmath}
\mu^{n+1}=  \frac{1}{2}\left ( \nu^n +\eta^n + \frac{1}{s}\left ( p^{n+1} - b^{n+1} \right )\right ) 
\end{dmath}
\begin{algorithm*}[!t]
	\begin{algorithmic}
	\State \textbf{Initialization}: $\left (\phi^n_0 , \nu^n_0 \right ) \in V \times H$ given arbitrarily.
	\Repeat
		\State (1) Compute $ \phi_k^n $ such that: $ \mathcal{L}_{p^n, \mu^n}\left ( \phi_k^n, \nu_{k-1}^n \right ) \leq  \mathcal{L}_{p^n, \mu^n}\left ( \phi, \nu_{k-1}^n \right ), \quad \forall \phi \in V $; \\
		\State (2) $ \nu^n_k= \nu_{k-1}^n+ \rho_{\nu} \left ( B \phi_{k}^n - p^n - s \left ( \nu_{k-1}^n - \mu^n \right ) \right ) $;
	\Until{Convergence}
	\end{algorithmic}
	\caption{For solving Step 1 in Algorithm \ref{alg:1}}
	\label{alg:1-1}
\end{algorithm*}

\begin{algorithm*}[!t]
	\begin{algorithmic}
	\State \textbf{Initialization}: $\left (q^n_0 , \eta^n_0 \right ) \in H \times H$ given arbitrarily.
	\Repeat
		\State (1) Compute $ q^n_k $ such that: $ \mathcal{L}_{ b^n, \mu^n}\left ( q^n_k, \eta^n_{k-1} \right ) \leq \mathcal{L}_{ b^n, \mu^n}\left ( q, \eta^n_{k-1} \right ), \quad \forall q \in H $;\\
		\State (2) $ \eta^n_k= \eta^n_{k-1}+ \rho_{\eta} \left (  b^n - q^n_{k} - s \left ( \eta^n_{k-1} - \mu^n \right ) \right ) $;
	\Until{convergence}
	\end{algorithmic}
	\caption{For solving Step 3 in Algorithm \ref{alg:1}}
	\label{alg:1-2}
\end{algorithm*}

The conditions on the problem data and the values of the parameters $\rho$, $\rho_\nu$ and $\rho_\eta$ for which Algorithm \ref{alg:1} and its sub-iterations converge are stated in the two following theorems. Here, the symbol $\to$ denotes strong convergence, whereas $\rightharpoonup$ denotes weak convergence.
\begin{theorem}[Convergence of Algorithm \ref{alg:1}]
	\label{thm:conv-alg1}
	Assume that
	\begin{enumerate}
		\item $G$ and $F$ are convex, lower semi-continuous functionals and uniformly convex on the bounded subsets of $V$ and $H$ respectively,
		\item $I$ is a proper, convex, lower semi-continuous functional on $H$,
		\item $B$ is a linear map from $V$ to $H$, that is injective and has a closed range,
		\item $\mathcal{L}_{r,s}$ admits a saddle point $\left ( \phi^*, q^*, p^*, b^*, \mu^*, \nu^*, \eta^* \right )$,
		\item $\mathcal{L}_{r,s}$ is coercive for $\left( \phi, q\right )$ with all other variables fixed, proper for $\phi$ with all other variables fixed and proper for $q$ with all other variables fixed.
	\end{enumerate}
	If $\rho$ verifies
	\begin{equation}
	0 < \rho < \frac{2rs^2 +s}{ \left ( 1+ rs \right )^2},
	\end{equation}
	then  $F\left ( q^n \right )+ G\left ( \phi^n \right )-  I\left ( \mu^n \right )+\left \langle \mu^n, p^n - b^n \right \rangle \to F\left ( q^* \right )+ G\left ( \phi^* \right ) -  I\left ( \mu^* \right )+\left \langle \mu^*, B\phi^*-  q^* \right \rangle$, $\phi^n \to \phi^*$, $p^n \to B \phi^*$, $q^n \to q^*$, $b^n \to q^*$, $\mu^n \rightharpoonup  \hat{\mu}$, $\nu^n \rightharpoonup  \hat{\mu}$ and $\eta^n \rightharpoonup  \hat{\mu}$, where $\left ( \phi^*, q^*, p^*, b^*, \hat{\mu}, \hat{\mu}, \hat{\mu} \right )$ is a saddle point of $\mathcal{L}_{r,s}$.
	
	Furthermore, if $I$ is uniformly convex on the subsets of $H$, then we have: $\mu^n \to  \mu^*$, $\nu^n \to  \mu^*$ and $\eta^n \to  \mu^*$.
\end{theorem}
\begin{proof}
	in supporting material.
\end{proof}
\begin{lemma}[Convergence of Algorithm \ref{alg:1-1}]
	\label{lem:conv-alg1-1}
	Under the same assumptions of Theorem \ref{thm:conv-alg1}, if $\rho_\nu$ verifies
	\begin{equation}
	0< \rho_{\nu} < \frac{2r}{2rs+1},
	\end{equation}
	then $\phi^n_k \to \phi^n$ and $\nu^n_k \to \nu^n$, where $\left( \phi^n, \nu^n \right )$ is a saddle point of $\mathcal{L}_{p^n, \mu^n}$.
\end{lemma}
\begin{proof}
	in supporting material.
\end{proof}
The convergence result for Algorithm \ref{alg:1-2} as well as its proof is mutatis mutandis the same as that of Algorithm \ref{alg:1-1}.

Although Algorithm \ref{alg:1} is guaranteed to converge, the cost of its inner steps 1 and 3 degrades its performance if the computation of $\phi^n$ and $q^n$ is not significantly less complex than that of $\mu^n$.
For most problems, however, Algorithm \ref{alg:2} is more likely to deliver a better performance.
In this algorithm, solving the two saddle point sub-problems is replaced by a minimization over $\phi$ and $q$ and a gradient ascent update rule for $\nu$ and $\eta$.   
Notice that in Step 2, the computation of $p^n$ requires the value $\nu^{n-1}$ and in turn, the computation of $\nu^n$  in Step 3 requires the value $p^{n-1}$. The same applies for Step 5 and 6.
This interdependence is necessary in the proof of the convergence result stated in Theorem \ref{thm:conv-alg2}, but the algorithm does convergence if the most recent values of the variables are used albeit for smaller step sizes.

\begin{algorithm*}[!t]
	\begin{algorithmic}
		\State \textbf{Initialization}: $\left ( p^0, b^0, \mu^0,\nu^0,\eta^0 \right ) \in H^5 $ given arbitrarily.
		\Repeat
			\State (1) Compute $ \phi^n $ such that: $ \mathcal{L}_{p^{n-1},\mu^n} \left ( \phi^n, \nu^{n-1} \right ) \leq \mathcal{L}_{p^{n-1},\mu^n} \left ( \phi, \nu^{n-1} \right ) $; \\ 
			\State (2)  $ p^n = p^{n-1} - \rho_r \left ( \mu^n - \nu^{n-1} + r \left ( p^{n-1} - B\phi^n \right) \right ) $;\\
			\State (3) $ \nu^n= \nu^{n-1} + \rho_s \left ( B \phi^n -p^{n-1} - s \left ( \nu^{n-1}- \mu^n \right) \right ) $;\\
			\State (4) Compute $ q^n $ such that: $ \mathcal{L}_{b^{n-1},\mu^n} \left ( q^n, \eta^{n-1} \right ) \leq \mathcal{L}_{b^{n-1},\mu^n} \left ( q, \eta^{n-1} \right ),\quad \forall q \in H $;\\
			\State (5) $ b^n = b^{n-1} - \rho_r \left ( \eta^{n-1} - \mu^n +r \left ( b^{n-1} - q^n \right ) \right ) $;\\
			\State (6) $ \eta^n= \eta^{n-1} + \rho_s \left ( b^{n-1} -q^n - s \left ( \eta^{n-1}- \mu^n \right) \right ) $;\\
			\State (7) $ \mu^{n+1}= \arg\inf_{\mu} [ s | \mu - ( \nu^n + \eta^n + ( p^n - b^n  ) / s ) / 2 |^2 + I ( \mu  ) ] $;
		\Until{convergence}
	\end{algorithmic}
	\caption{For solving Problem \eqref{eq:aug-lagrangian}}
	\label{alg:2}
\end{algorithm*}

Theorem \ref{thm:conv-alg2} states the conditions on $\rho_r$ and $\rho_s$ under which Algorithm \ref{alg:2} converges.
\begin{theorem}[Convergence of Algorithm \ref{alg:2}]
	\label{thm:conv-alg2}
	Under the same assumptions of Theorem \ref{thm:conv-alg1}, if $\rho_r$, $\rho_s>0$ satisfy
	\begin{dmath}
	\label{eq:alg2-system}
	\left\{\begin{array}{lr}
	2  s - \rho_r - \rho_s s^2 - \left | \rho_r r - \rho_s s  \right | & >0, \\
	2 r - \rho_r r^2 - \rho_s - \left | \rho_r r - \rho_s s  \right | & >0,
	\end{array}\right.
	\end{dmath}
	then the same convergence results of Theorem \ref{thm:conv-alg1} are obtained.
\end{theorem}
\begin{proof}
	in supporting material.
\end{proof}
It is easy to check that System \eqref{eq:alg2-system} always admits a solution for any strictly positive values of $r$ and $s$.

Numerical experiments show that it is possible to speedup the convergence of Algorithm \ref{alg:2} in term of number of iterations by updating $\mu$ two times per iteration instead of one as shown in Algorithm \ref{alg:3}.
However, this should be weighted against the cost increase per iteration resulting from solving an additional minimization problem.

\begin{algorithm*}[!t]
	\begin{algorithmic}
		\State \textbf{Initialization}: $\left ( p^0, b^0, \mu^0,\nu^0,\eta^0 \right ) \in H^5$ given arbitrarily.
		\Repeat
				\State (1) Compute $ \phi^n $ such that: $ \mathcal{L}_{p^{n-1},\mu^n} \left ( \phi^n, \nu^{n-1} \right ) \leq \mathcal{L}_{p^{n-1},\mu^n} \left ( \phi, \nu^{n-1} \right ), \quad \forall \phi \in V $;\\
				\State (2)  $ p^n = p^{n-1} - \rho_r \left ( \mu^n - \nu^{n-1} + r \left ( p^{n-1} - B\phi^n \right) \right ) $; \\
				\State (3) $ \nu^n= \nu^{n-1} + \rho_s \left ( B \phi^n -p^{n-1} - s \left ( \nu^{n-1}- \mu^n \right) \right ) $; \\
				\State (4) $ \mu^{n+1/2}= \arg\inf_{\mu \in H} [  s  | \mu -  ( \nu^n +\eta^{n-1} +  ( p^n - b^{n-1}  )/s )/2  |^2 + I ( \mu  ) ] $; \\
				\State (5) Compute $ q^n $ such that: $ \mathcal{L}_{b^{n-1},\mu^{n+1/2}} \left ( q^n, \eta^{n-1} \right ) \leq \mathcal{L}_{b^{n-1},\mu^{n+1/2}} \left ( q, \eta^{n-1} \right ),\quad \forall q \in H $; \\
 				\State (6) $ b^n = b^{n-1} - \rho_r \left ( \eta^{n-1} - \mu^{n+1/2} +r \left ( b^{n-1} - q^n \right ) \right ) $;\\
				\State (7) $ \eta^n= \eta^{n-1} + \rho_s \left ( b^{n-1} -q^n - s \left ( \eta^{n-1}- \mu^{n+1/2} \right) \right ) $;\\
				\State (8) $ \mu^{n+1}= \arg\inf_{\mu \in H} [ s  | \mu - ( \nu^n +\eta^n +  ( p^n - b^n  )/s  )/2  |^2 + I ( \mu ) ] $.
		\Until{convergence}
	\end{algorithmic}
	\caption{For solving Problem \eqref{eq:aug-lagrangian}}
	\label{alg:3}
\end{algorithm*}

As stated earlier, Problem \eqref{eq:original-lagrangian} is a special case of Problem \eqref{eq:general-lagrangian}, but, from a theoretical standpoint, it lacks coercivity and uniform convexity.
This can be dealt with by adding small quadratic terms in $q$ and $\phi$.
In practice, however,  the algorithms presented solve Problem \eqref{eq:original-lagrangian} successfully without such perturbation.
Using a similar argument to that used in  \cite{benamou00cfd}, computing $\phi^n$ can be shown to be  equivalent to solving the Poisson equation 
\begin{dmath}
-r \Delta \phi= \nabla \cdot \left ( \nu^{n-1} - r p^{n-1} \right ),
\end{dmath}
with Neumann boundary conditions in the time dimension given by 
\begin{dmath}
r\partial_t \phi \left ( 0, \cdot \right ) = \rho_0 - \nu_0^{n-1}\left ( 0, \cdot \right ) + r p^{n-1}_0\left ( 0, \cdot \right )
\end{dmath}
and 
\begin{dmath}
r\partial_t \phi \left ( 1, \cdot \right ) = \rho_1 - \nu_0^{n-1}\left ( 1, \cdot \right ) + r p^{n-1}_0\left ( 1, \cdot \right ),
\end{dmath}
where $\nu_0$ and  $p_0$  are the first elements of the vectors $\nu$ and $p$ respectively. 
Similarly, finding $q^n$ amounts to solving: 
\begin{dmath}
\label{eq:compute-q-optimal-transport-alg2}
\inf_{q \in K} \left | b^{n-1} + \frac{\eta^{n-1}}{r} - q \right |^2.
\end{dmath}
The expressions for $\phi^n$ and $q^n$ given here correspond to Algorithm \ref{alg:2}. The translation to Algorithm 1 and Algorithm \ref{alg:3} is straightforward and requires only the appropriate change of indices.

\section{Numerical Experiments}
\label{sec:experiments}

To illustrate the working of the proposed algorithms, a number of experiments are presented in this section. Namely, we will use the proposed algorithms to solve several example problems. We will also compare the time performance of the algorithms on a sample problem and investigate the effects of the parameters $ r $ and $ s $ on convergence. 

In all these experiments, and unless otherwise stated, $r$ and $s$ are chosen to be 1, for Algorithm \ref{alg:1}: $\rho= \rho_\nu= \rho_\eta= 0.5$, and for Algorithm \ref{alg:2} and \ref{alg:3}: $\rho_r=0.4$ $\rho_s=1$.

Before presenting the experimental results, a brief discussion of the implementation of the proposed algorithms is first presented.

\subsection{Implementation}
In the examples presented in this section, the space domains are two-dimensional and are discretized using a regular square grid.  A seven-point stencil is used to discretize the time-space Laplacian operator used in the Poisson equation. The linear system resulting from this discretization is solved using the BoomerAMG solver \cite{HensonVanEmden2002BApa}, which is a parallel solver part of the hypre library \cite{FalgoutRobertD2006Ciih}. The optimization problem for computing $ q $ is solved analytically, and the resulting Karush-Kuhn-Tucker optimality equations are solved using the Jenkins-Traub algorithm \cite{Jenkins70}.

The proposed algorithms are implemented in C++ and use distributed parallelism with MPI for communication and synchronization. This allows to accelerate computations and to take full advantage of the capabilities offered by the hypre library.

\subsection{A transport problem with bound constraint}
\begin{figure*}[!t]
	\begin{center}
		\includegraphics[width=12cm]{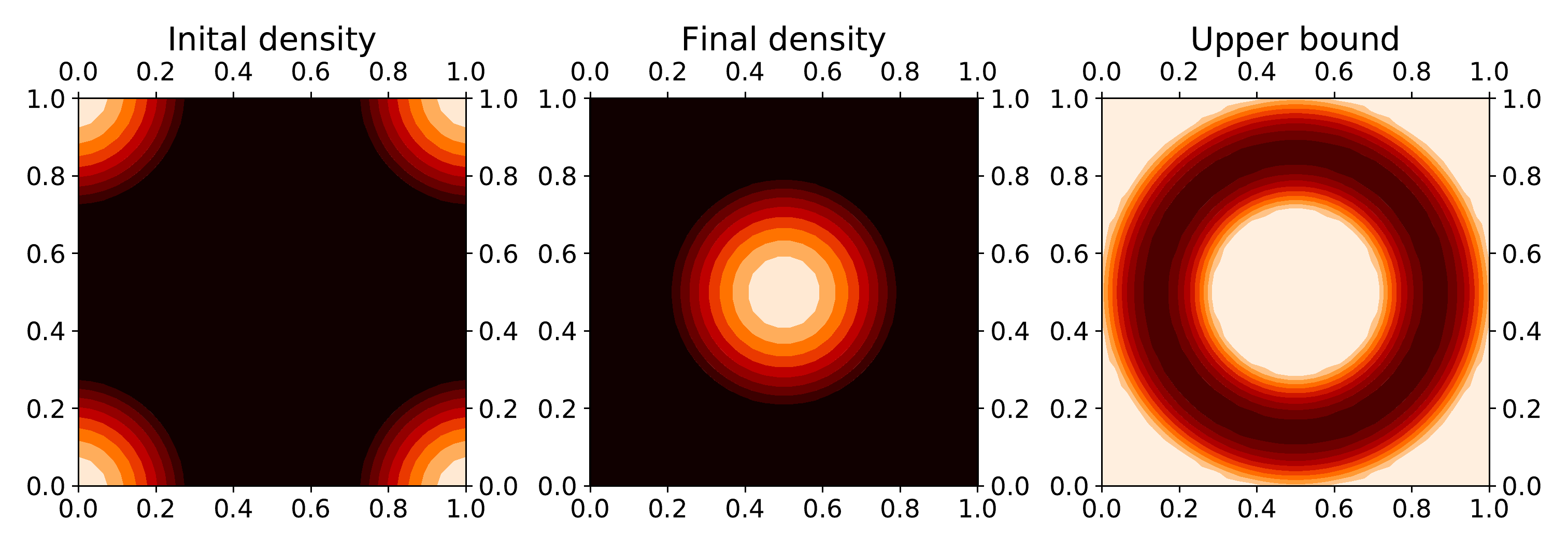}
	\end{center}
	\caption{This problem consists in transporting the initial density (shown on the left) to match the final density (shown in the middle) under the constraint that the density does not exceed the upper bound (shown on the top right).}
	\label{fig:bounded}
\end{figure*}

In this problem, each point in the space has an upper bound on its density, which may vary from one point to another, and the goal is transport the initial density to the final one without violating any bound constraint.  As shown in Fig. \ref{fig:bounded}, the imposed upper bound creates a barrier around the inner area of the domain preventing the mass from flowing freely towards the center. The space domain is discretized using a $ 64\times64 $ square grid with periodic boundaries, whereas the time resolution is 65. Notice that, since the mass is required to move to the center of the domain, the periodicity of the domain has little to no effect on the resulting flow.

\begin{figure*}[!t]
	\begin{center}
		\includegraphics[width=16cm]{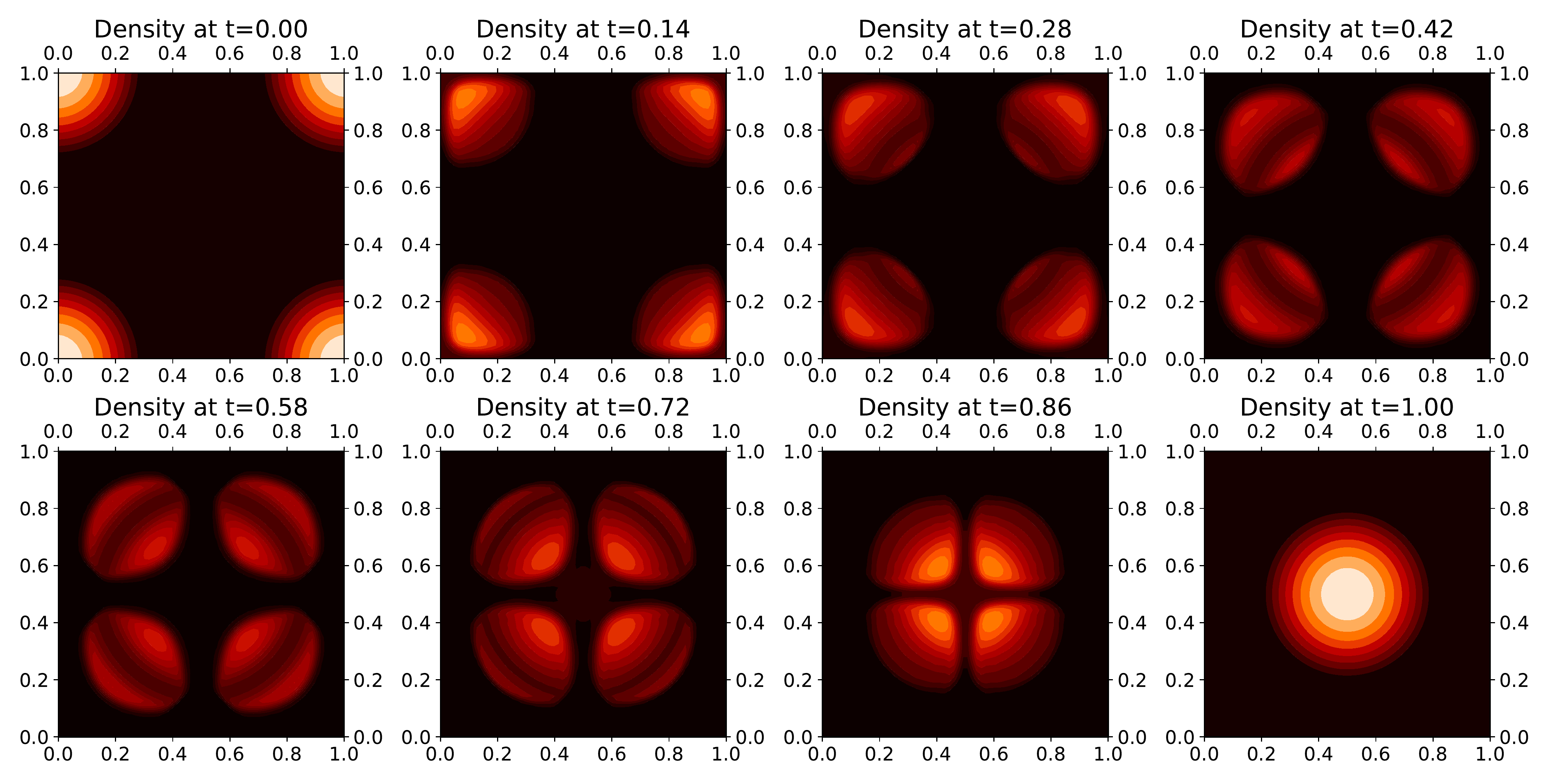}
	\end{center}
	\caption{The solution to the constrained transportation problem shown in Fig. \ref{fig:bounded}.}
	\label{fig:bounded-sol}
\end{figure*}

\begin{figure*}[!t]
	\begin{center}
		\includegraphics[width=16cm]{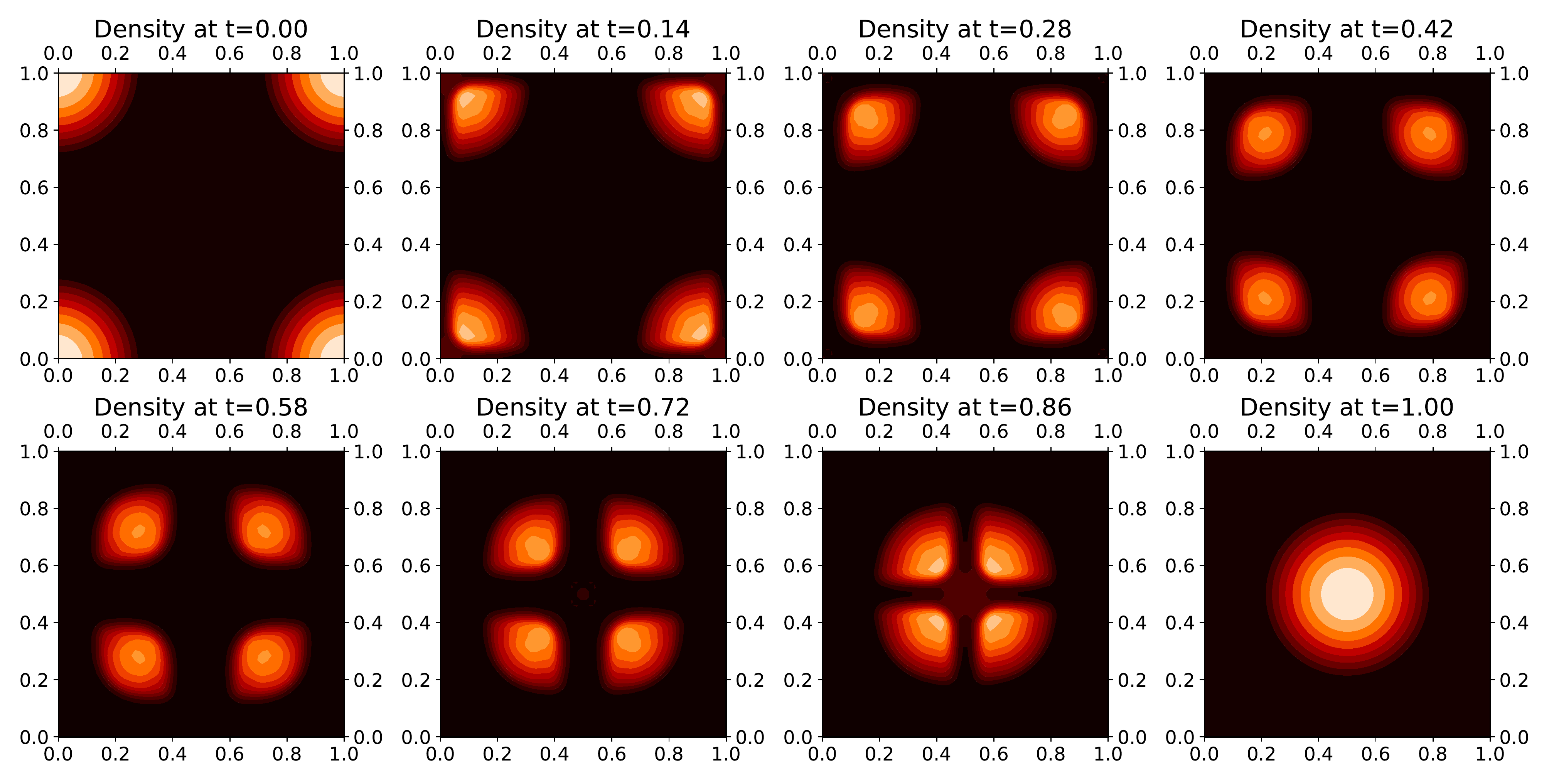}
	\end{center}
	\caption{The solution to the transportation problem shown in Fig. \ref{fig:bounded} without considering the constraint (that is, without imposing an upper bound on the density).}
	\label{fig:nobound-sol}
\end{figure*}

\begin{figure*}[!t]
	\begin{center}
		\includegraphics[width=16cm]{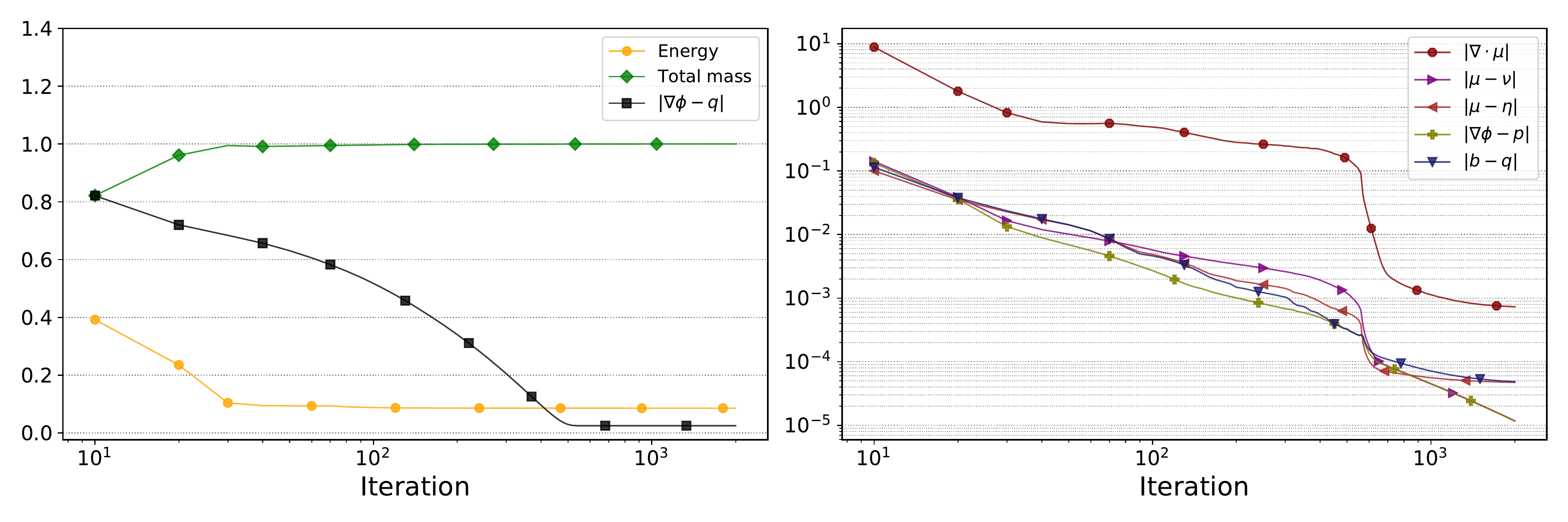}
	\end{center}
	\caption{A log-log plot of the evolution of the convergence criteria for the transport problem shown in Fig. \ref{fig:bounded} when solved by Algorithm \ref{alg:3}. The quantities plotted on the left side converge to non-zero values, whereas those on the right converge to zero.}
	\label{fig:bound-conv}
\end{figure*}

Notice in the solution shown in Fig. \ref{fig:bounded-sol} how the mass is ''squeezed'' under the barrier imposed by the upper bound before reaching its final destination in the center of the domain. For comparison, Fig. \ref{fig:nobound-sol} shows the optimal transport solution without the bound constraint; The mass is simply moved diagonally towards the center of the domain. 
Fig. \ref{fig:bound-conv} shows the evolution of the convergence criteria when running Algorithm \ref{alg:3} on this problem.
As expected, because of the constraint, $\left | \nabla \phi -q  \right | $ does not converge to zero. 

Bounds on the density can be useful in many application domains where the system cannot undertake severe compression. For instance, in a crowd evacuation application, imposing an upper bound on the density is essential to ensure the safety of crowd during the evacuation process. 

\subsection{A problem with a penalty on momentum}
\begin{figure*}[!t]
	\begin{center}
		\includegraphics[width=12cm]{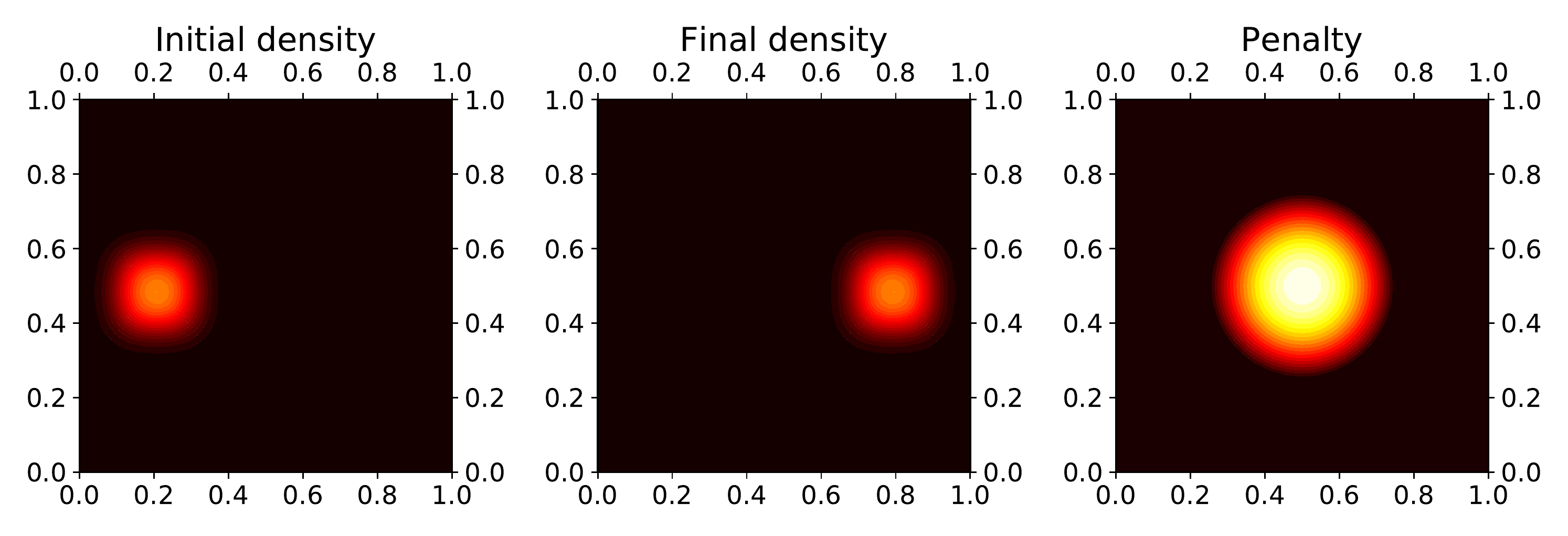}
		
		\vspace{0.2cm}
		
		\includegraphics[width=16cm]{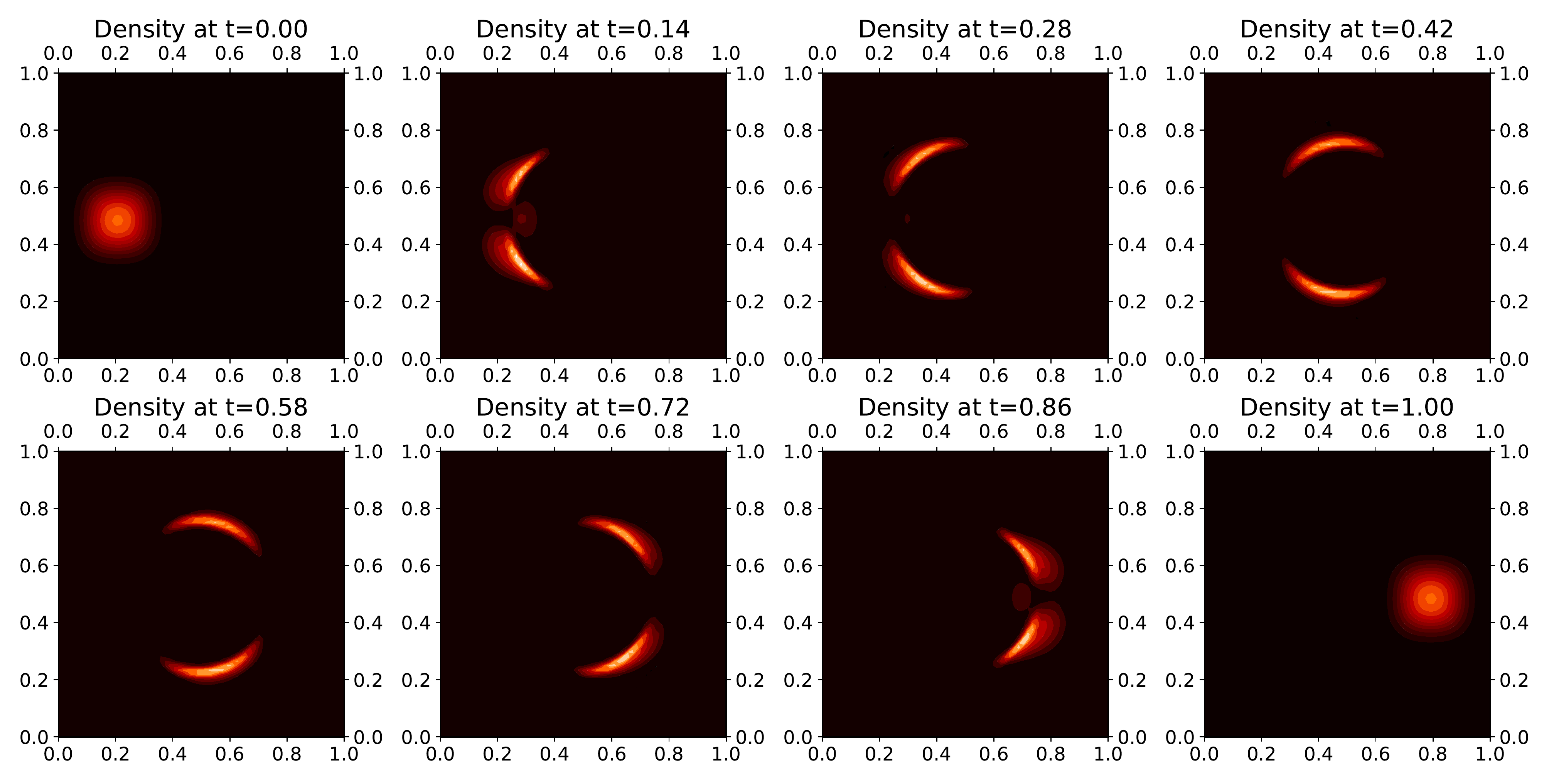}
	\end{center}
	\caption{This problem consists in transporting the initial density (shown on the top left) to match the final density (shown in the top middle) with a penalty imposed on transporting mass through the center of the domain (see top right).}
	\label{fig:penalty}
\end{figure*}

In this problem, a penalty is imposed on the momentum field at the center of the domain. This penalty is of the form $ \int_{0}^{1}\int_D \psi(x) \left|m(t,x)\right|^2 dxdt$, where $ \psi(x) $ is a non-negative continuous function which is null everywhere except within a disk area at center of the domain (see Fig. \ref{fig:penalty}). Although the density is not constrained, the penalty on the momentum field will cause the flow to avoid the central area. For this problem, we use a $ 65\times 65 $ square grid for discretizing the space domain and a time resolution of 65.  We impose a null-Neumann boundary condition on the space domains, which simply means that the domain is isolated and mass cannot flow through its boundary.
As the solution shows, the mass is essentially transported around the central area to avoid penalty.

More complex penalty functions on the momentum can be imposed to capture various constraints imposed by the physical world. For instance, one can impose a given direction fo the flow at different parts of the domain, or impose a limit on the norm of the momentum, which may reflect physical limitations in the system. Notice that the non-smoothness of the constraint translates into a non-smooth flow of the density.

\subsection{A problem with a area of constant density}
In this experiment, a constrained problem where the mass density must remain constant inside a disk area at the center of the domain is solved. Here, we also use  a $ 65\times 65 $ square grid for discretizing the space domain, impose a null-Neumann boundary condition on space, and use a time resolution of 65. 

\begin{figure*}[!t]
	\begin{center}
		\includegraphics[width=16cm]{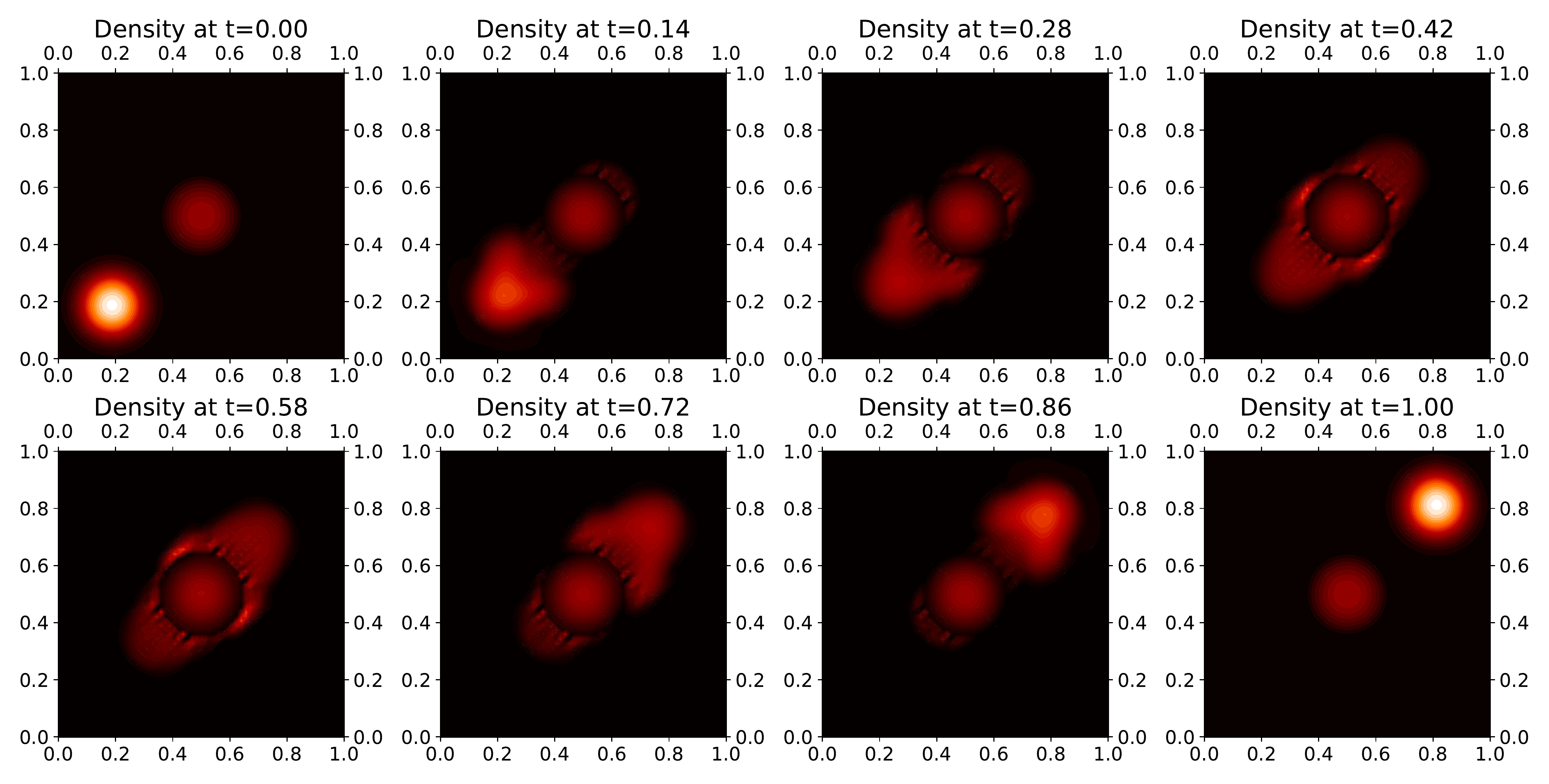}
	\end{center}
	\caption{This problem consists in transporting the initial density ($ t=0 $) to match the final density ($ t=1 $) under the constraint that the density at the center of the domain remains unchanged during the transportation process.}
	\label{fig:const}
\end{figure*}

Fig. \ref{fig:const} shows the solution computed using Algorithm \ref{alg:2}. The mass is transported through the center of the domain without changing its density, while the extra mass is routed around the constraint area. 

\subsection{Time performance comparison}
The next experiment gives a performance comparison between the three algorithms in terms of computation time. For this, we run the three algorithms on a sample problem, namely the constrained momentum problem shown in Fig. \ref{fig:penalty}, and report in Table \ref{tab:time} the time and number of iterations taken to solve the problem. The stopping criterion is that the change in density between two consecutive iterations is smaller than a certain threshold. More precisely, we stop the algorithms at iteration $ k $ where $\left|\rho^{k}-\rho^{k-1}\right| < \epsilon =1.0e-3 $. The time reported is the average of 10 runs, where each run consists of a parallel execution of the algorithm on four cores of an Intel Core i7-8750H CPU @ 2.20GHz. The space domain is discretized as a regular grid of size 33$\times$33 and the time resolution is 33. 

The results show that for this problem, Algorithm \ref{alg:3} reaches convergence in almost half the time and  number of iterations required by Algorithm \ref{alg:2} thanks to the double update of $ \mu $. Regarding Algorithm \ref{alg:1}, although it takes the same number of iterations as Algorithm \ref{alg:2}, the time it requires is much longer because of the time spent in repeatedly solving the linear system associated with the Poisson equation required for computing $ \phi $. 

Notice that, for this problem, the time spent in solving this linear system dominates the computation time of the algorithms, and in particular, relatively to that of computing $ \mu $. This explains the superiority of Algorithm \ref{alg:3} over Algorithm \ref{alg:2}. However, in problems where the optimization problem associated with $ \mu $ is more computationally intensive than that of finding $ \phi $, the situation might be reversed.

\begin{table}[!t]
	\centering
	\caption{Time, with relative standard deviation, and number of iterations taken by Algorithm \ref{alg:1}, \ref{alg:2} and \ref{alg:3} to solve the transport problem shown in Fig. \ref{fig:penalty}.}
	\label{tab:time}
	\begin{tabular}{lccc}
		\toprule
		& Algorithm \ref{alg:1}& Algorithm \ref{alg:2}& Algorithm \ref{alg:3} \\
		\midrule
		Time (sec) &1248.22 $\pm$3\% & 169.25 $\pm$5\% &  87.78 $\pm$9\%\\
		Number of iterations &841&841&450\\
		\bottomrule
	\end{tabular}
\end{table}

\subsection{Effect of the parameters $ r $ and $ s $}
\begin{figure*}[!t]
	\begin{center}
		\includegraphics[width=12cm]{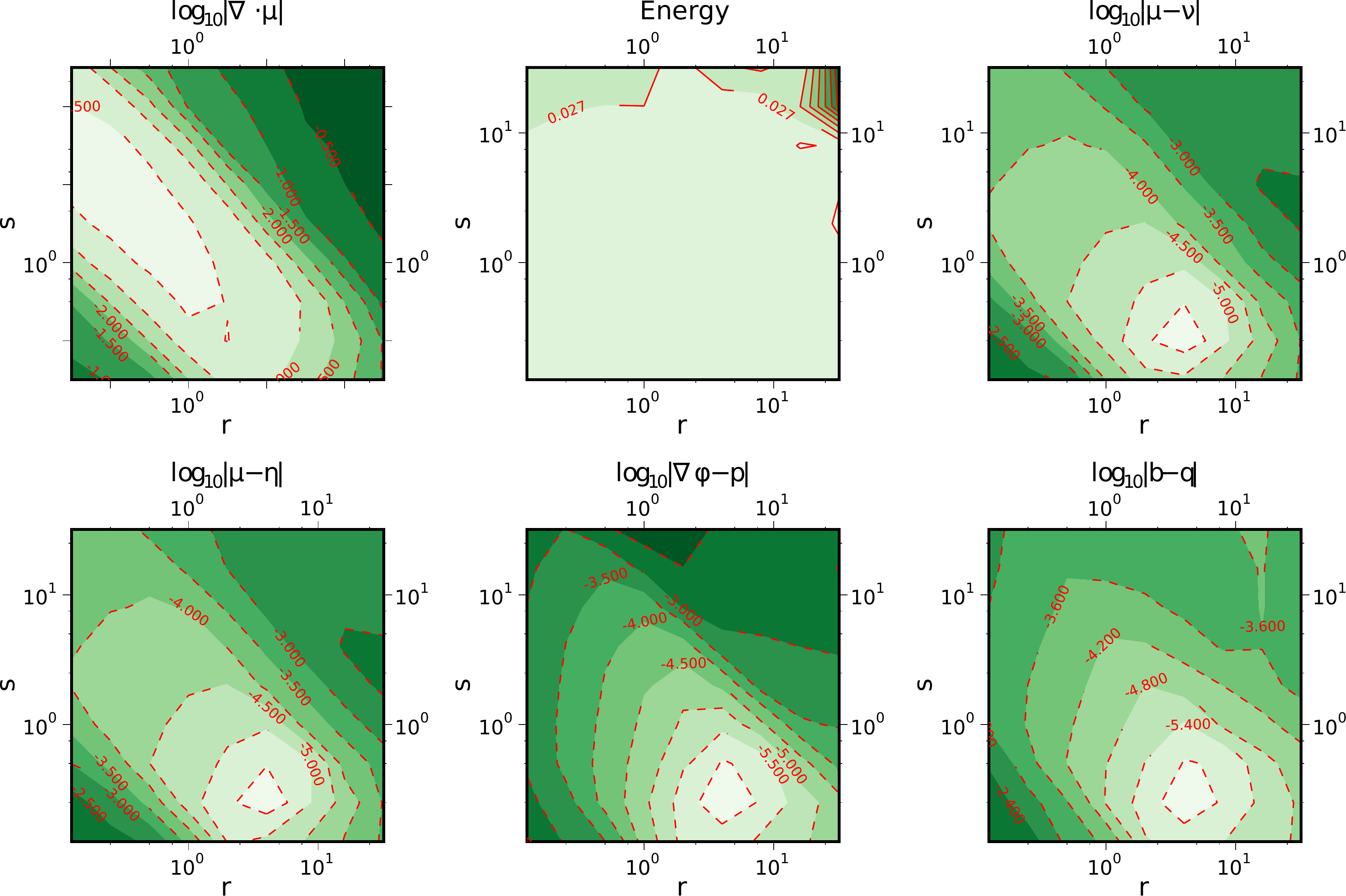}
	\end{center}
	\caption{Effect of the values of the parameters $r$ and $s$ on the convergence of Algorithm \ref{alg:2}. The plotted results are obtained after 500 iterations. The plots are all logarithmic in $r$ and $s$. }
	\label{fig:rs}
\end{figure*}

To investigate the effect of the parameters $ r $ and $ s $ on convergence, we plot the different convergence criteria obtained after 500 iterations of Algorithm \ref{alg:2} on the bound constraint problem for different values of $ r $ and $ s $. As shown in Fig. \ref{fig:rs}, except for extreme values, $ r $ and $ s $ have little effect of the energy. The remaining criteria, on the other hand, exhibit clear sensitivity to changes in $ r $ and $ s $ with $ \left|\nabla \cdot \mu\right| $ favoring small values of $ r $ and large values of $ s $, and vice versa for the other criteria.

\section{Conclusion}

From the theoretical perspective, there are several issues related to constrained optimal transport that deserve further investigation, for instance, the existence of minimizers and the regularity of the solutions. We believe that this endeavor can be best tackled using the framework of Lagrangian action introduce in  \cite{villani09}.  It is, however, worth mentioning that, even for free optimal transport, regularity is still an active area of research with many unsettled problems \cite{villani09}.

On the computational side, the requirement by the algorithms to solve an optimization problem for finding $ \mu $ at each iteration can be time consuming, especially for problems with complex constraints such those defined by a large number of inequalities. A possible solution is to replace the minimization problem by a gradient descent update rule for the case of a soft constraint $ I $ that is Frechet differentiable and an interior point technique when $ I $ is hard. 

From the application stand point, and in addition to the classical application in image registration, the possibility of constraining the transport process can be used to probe the intrinsic geometry of the feasible set. This allows to compute the intrinsic distance between probability distributions, which may greatly differ from that in the ambient space. Since the concept of distance, or dissimilarity, between data plays a key role in virtually all areas of machine learning and pattern recognition, the flexibility in designing new distances by considering different constraints can play an important role in improving the performance of existing learning algorithms.


%

\appendices

\ifCLASSOPTIONcompsoc
  \section*{Acknowledgments}
\else
  \section*{Acknowledgment}
\fi

\ifCLASSOPTIONcaptionsoff
  \newpage
\fi



\bibliographystyle{IEEEtran}
%

%




\includepdf[pages=-]{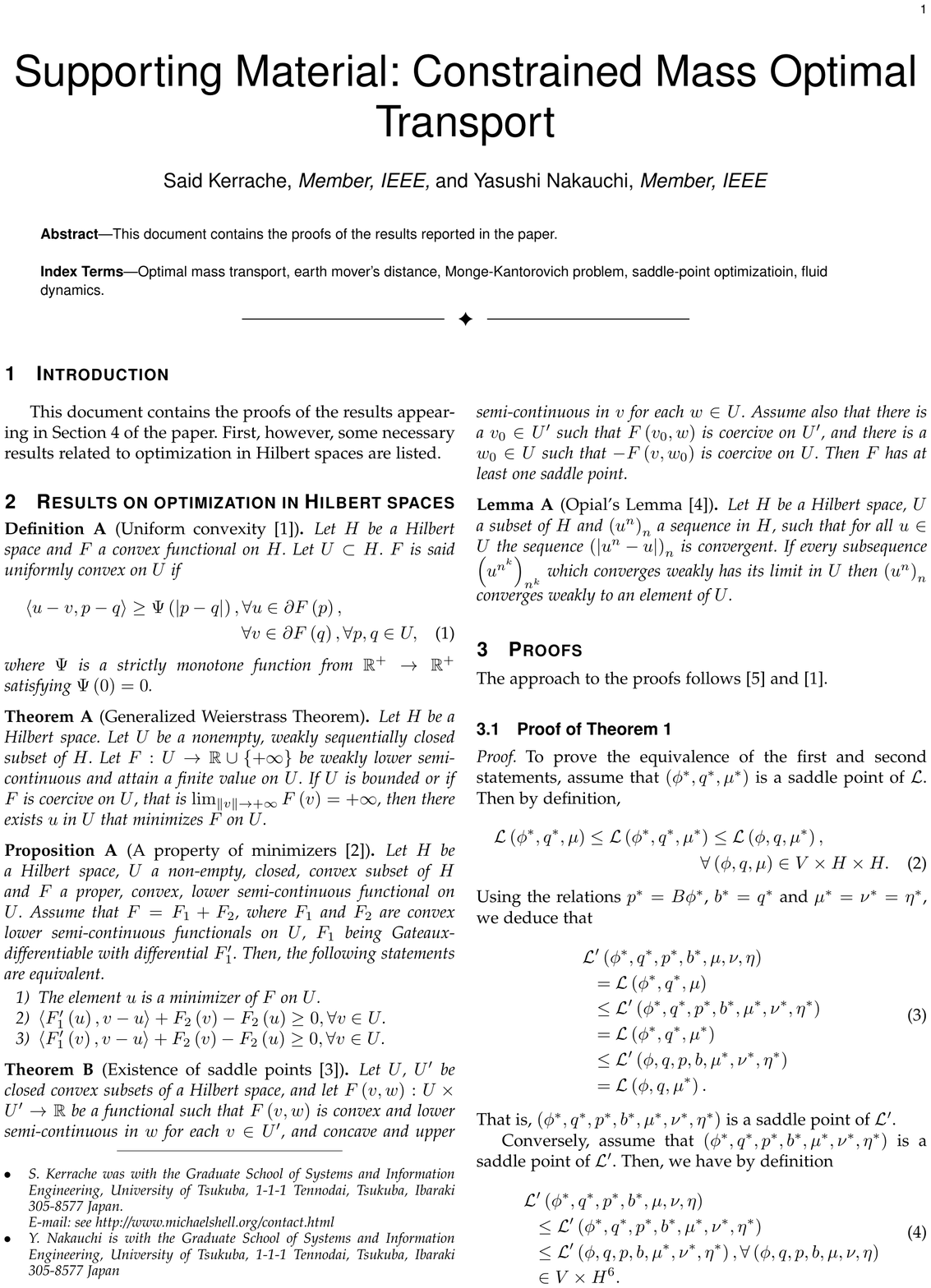}
\end{document}


\title{Supporting Material: Constrained Mass Optimal Transport}

\author{Said~Kerrache,~\IEEEmembership{Member,~IEEE,}
        and~Yasushi~Nakauchi,~\IEEEmembership{Member,~IEEE}
\IEEEcompsocitemizethanks{\IEEEcompsocthanksitem S. Kerrache was with the Graduate School of Systems and Information Engineering, University of Tsukuba, 1-1-1 Tennodai, Tsukuba, Ibaraki 305-8577 Japan.\protect\\
E-mail: see http://www.michaelshell.org/contact.html
\IEEEcompsocthanksitem Y. Nakauchi is with the Graduate School of Systems and Information Engineering, University of Tsukuba, 1-1-1 Tennodai, Tsukuba, Ibaraki 305-8577 Japan}
}

\markboth{}%
{}

\IEEEtitleabstractindextext{%
\begin{abstract}
This document contains the proofs of the results reported in the paper.
\end{abstract}

\begin{IEEEkeywords}
Optimal mass transport, earth mover's distance, Monge-Kantorovich problem, saddle-point optimizatioin, fluid dynamics. 
\end{IEEEkeywords}}

\maketitle

\IEEEdisplaynontitleabstractindextext

\IEEEpeerreviewmaketitle

\IEEEraisesectionheading{\section{Introduction}\label{sec:introduction}}

This document contains the proofs of the results appearing in Section 4 of the paper.
First, however, some necessary results related to optimization in Hilbert spaces are listed.

\section{Results on optimization in Hilbert spaces}

\begin{definition}[Uniform convexity \cite{glowinski89}]
	\label{def:uniformly-convex}
	Let $H$ be a Hilbert space and $F$ a convex functional on $H$.
	Let $U \subset H$.
	$F$ is said uniformly convex on $U$ if
	\begin{multline}
	\left \langle u- v, p-q \right \rangle \geq \Psi \left ( \left | p-q \right | \right ), \forall u \in \partial F\left ( p \right ), \\ \forall v \in \partial F\left ( q \right ) , \forall p,q \in U,
	\end{multline}
	where $\Psi$ is a strictly monotone function from $\R^+ \to \R^+$ satisfying $\Psi \left (0 \right ) =0$.
\end{definition}
\begin{theorem}[Generalized Weierstrass Theorem]
	\label{thm:weiestrass}
	Let $H$ be a Hilbert space.
	Let $U$ be a nonempty, weakly sequentially closed subset of $H$.
	Let $F: U \to \R \cup \left \{ +\infty \right \}$ be weakly lower semi-continuous and attain a finite value on $U$.
	If $U$ is bounded or if $F$ is coercive on $U$,  that is $\lim_{\left \| v \right \| \to +\infty} F\left ( v \right )= + \infty$, then there exists $u$ in $U$ that minimizes $F$ on $U$.
\end{theorem}
\begin{proposition}[A property of minimizers \cite{ekeland74convex}]
	
	Let $H$ be a Hilbert space, $U$ a non-empty, closed, convex subset of $H$ and $F$ a proper, convex, lower semi-continuous functional on $U$.
	Assume that $F= F_1 + F_2$, where $F_1$ and $F_2$ are convex lower semi-continuous functionals on $U$, $F_1$ being Gateaux-differentiable with differential $F_1'$.
	Then, the following statements are equivalent.
	\begin{enumerate}
		\item The element $u$ is a minimizer of $F$ on $U$.
		\item $\left \langle F_1'\left ( u \right ), v-u \right \rangle + F_2\left ( v \right ) - F_2\left ( u \right ) \geq 0, \forall v \in U$.
		\item $\left \langle F_1'\left ( v \right ), v-u \right \rangle + F_2\left ( v \right ) - F_2\left ( u \right ) \geq 0, \forall v \in U$.
	\end{enumerate}
	
\end{proposition}
\begin{theorem}[Existence of saddle points \cite{schechter05}]
	\label{thm:existence-saddle-point}
	Let $U$, $U'$ be closed convex subsets of a Hilbert space, and let $F\left( v, w \right ): U \times U' \to \R$ be a functional such that  $F \left( v, w \right) $ is convex and lower semi-continuous in $w$ for each $v \in U'$, and concave and upper semi-continuous in $v$ for each $w \in U$.
	Assume also that there is a $v_0 \in U'$  such that $F\left ( v_0,w  \right )$ is coercive on $U'$, and there is a $w_0 \in U$  such that $-F\left ( v,w_0  \right )$  is coercive on $U$.
	Then $F$ has at least one saddle point.
\end{theorem}
\begin{lemma}[Opial's Lemma \cite{glowinski81}]
	\label{lem:opial}
	Let $H$ be a Hilbert space, $U$ a subset of $H$ and $\left (u^n\right )_n$  a sequence in $H$, such that for all $u \in U$ the sequence $\left ( \left | u^n -u  \right | \right )_{n}$ is convergent.
	If every subsequence $\left ( u^{n^k} \right )_{n^k}$ which converges weakly has its limit in $U$ then $\left (u^n\right )_n$ converges weakly to an element of $U$.
\end{lemma}
\section{Proofs}
The approach to the proofs follows \cite{fortin83} and \cite{glowinski89}.

\subsection{Proof of Theorem 1}
\begin{proof}[Proof]
	To prove the equivalence of the first and second statements, assume that $\left ( \phi^*, q^*, \mu^* \right )$ is a saddle point of $\mathcal{L}$. Then by definition,
	\begin{multline}
	\mathcal{L} \left ( \phi^*, q^*, \mu \right ) \leq \mathcal{L} \left ( \phi^*, q^*, \mu^* \right ) \leq \mathcal{L} \left ( \phi, q, \mu^* \right ), \\ \forall \left( \phi, q, \mu \right) \in V \times H \times H.
	\end{multline}
	Using the relations $p^*= B \phi^*$, $b^*=q^*$ and $\mu^*=\nu^*=\eta^*$, we deduce that
	\begin{dmath}
		\mathcal{L}'\left ( \phi^*, q^*, p^*, b^*, \mu, \nu, \eta \right) = \mathcal{L}\left ( \phi^*, q^*, \mu \right)  \leq \mathcal{L}'\left ( \phi^*, q^*, p^*, b^*, \mu^*, \nu^*, \eta^* \right)	= \mathcal{L}\left ( \phi^*, q^*, \mu^* \right)  \leq \mathcal{L}'\left ( \phi, q, p, b, \mu^*, \nu^*, \eta^* \right) = \mathcal{L}\left ( \phi, q, \mu^* \right).
	\end{dmath}
	That is, $ \left ( \phi^*, q^*, p^*, b^*, \mu^*, \nu^*, \eta^* \right) $ is a saddle point of  $\mathcal{L}'$.
	
	Conversely, assume that  $\left ( \phi^*, q^*, p^*, b^*, \mu^*, \nu^*, \eta^* \right) $ is a saddle point of  $\mathcal{L}'$. Then, we have by definition
	\begin{dmath}
		\label{eq:sadlle-point-l-prime}
		\mathcal{L}'\left ( \phi^*, q^*, p^*, b^*, \mu, \nu, \eta \right) \leq \mathcal{L}'\left ( \phi^*, q^*, p^*, b^*, \mu^*, \nu^*, \eta^* \right) 
		\leq \mathcal{L}'\left ( \phi, q, p, b, \mu^*, \nu^*, \eta^* \right), \forall \left ( \phi, q, p, b, \mu, \nu, \eta \right) \in V \times H^6.
	\end{dmath}
	From the first inequality in \eqref{eq:sadlle-point-l-prime}, and by choosing $\mu= \mu^*$ and $\eta=\eta^*$, it follows that
	\begin{equation}
	\left \langle \nu , B\phi^* - p^* \right \rangle \leq \left \langle \nu^* , B\phi^* - p^* \right \rangle, \quad \forall \nu \in H,
	\end{equation}
	which implies that $B\phi^*=p^*$. Choosing this time $\mu= \mu^*$ and $\nu=\nu^*$, it follows that
	\begin{equation}
	\left \langle \eta , b^* - q^* \right \rangle \leq \left \langle \eta^* , b^* - q^* \right \rangle, \quad \forall \eta \in H,
	\end{equation} 
	yielding $b^*=q^*$.
	A similar argument using the second inequality in \eqref{eq:sadlle-point-l-prime} gives us $\mu^*=\nu^*=\eta^*$.
	The first equality results from choosing $q=q^*$, $\phi=\phi^*$ and $b=b^*$, whereas the second equality follows from the choice $q=q^*$, $\phi=\phi^*$ and $p=p^*$.
	These equalities imply
	\begin{dmath}
		\mathcal{L}\left ( \phi^*, q^*, \mu \right) = \mathcal{L}'\left ( \phi^*, q^*, p^*, b^*, \mu, \nu, \eta \right) 
		\\ \leq \mathcal{L}'\left ( \phi^*, q^*, p^*, b^*, \mu^*, \nu^*, \eta^* \right)
		= \mathcal{L}\left ( \phi^*, q^*, \mu^* \right)  \leq \mathcal{L}\left ( \phi, q, \mu^* \right)=  \mathcal{L}'\left ( \phi, q, p, b, \mu^*, \nu^*, \eta^* \right),
		\forall \left ( \phi, q, p, b, \mu, \nu, \eta \right) \in V \times H^6.
	\end{dmath}
	It follows that $\left ( \phi^*, q^*, \mu^* \right )$ is a saddle point of $\mathcal{L}$.
	This completes the proof of the first equivalence.
	
	Now we shall prove the equivalence of the second and third statements.
	Assume that $\left ( \phi^*, q^*, p^*, b^*, \mu^*, \nu^*, \eta^* \right) $ is a saddle point of  $\mathcal{L}'$. Then
	\begin{dmath}
		\mathcal{L}_{r,s} \left ( \phi^*, q^*, p^*, b^*, \mu, \nu, \eta \right ) \leq \mathcal{L}' \left ( \phi^*, q^*, p^*, b^*, \mu, \nu, \eta \right )  
		\leq \mathcal{L}' \left ( \phi^*, q^*, p^*, b^*, \mu^*, \nu^*, \eta^* \right  )   = \mathcal{L}_{r,s} \left ( \phi^*, q^*, p^*, b^*, \mu^*, \nu^*, \eta^* \right ) 
		\leq \mathcal{L}' \left ( \phi, q, p, b, \mu^*, \nu^*, \eta^* \right ) \\ \leq \mathcal{L}_{r,s} \left ( \phi, q, p, b, \mu^*, \nu^*, \eta^* \right )
		\forall \left ( \phi, q, p, b, \mu, \nu, \eta \right) \in V \times H^6. 
	\end{dmath}
	That is, $\left ( \phi^*,q^*, p^*, b^*, \mu^*, \nu^*, \eta^* \right )$ is also a saddle point of $\mathcal{L} _{r,s}$.
	
	Conversely, assume that $\left ( \phi^*,q^*, p^*, b^*, \mu^*, \nu^*, \eta^* \right )$ is a saddle point of $\mathcal{L} _{r,s}$. Then, by definition,
	\begin{dmath}
		\mathcal{L}_{r,s} \left ( \phi^*, q^*, p^*, b^*, \mu, \nu, \eta \right ) \leq \mathcal{L}_{r,s} \left ( \phi^*, q^*, p^*, b^*, \mu^*, \nu^*, \eta^* \right ) 
		\leq \mathcal{L}_{r,s} \left ( \phi, q, p, b, \mu^*, \nu^*, \eta^* \right ),  \forall \left ( \phi, q, p, b, \mu, \nu, \eta \right) \in V \times H^6.
	\end{dmath}
	The optimality conditions for $\nu^*$, $p^*$, $\eta^*$ and $b^*$ read as follows:
	\begin{equation}
	\label{eq:optimality-l-rs}
	\left\{\begin{array}{ll}
	\nu^*: & B \phi^* - p^* = s \left ( \nu^* - \mu^* \right ), \\
	p^*:  & r \left (B \phi^* - p^* \right )= \mu^* - \nu^*, \\
	\eta^*: & b^* - q^* = s \left ( \eta^* - \mu^* \right ), \\
	b^*: & r \left (b^* - q^* \right ) = \mu^* - \eta^*.
	\end{array} \right.
	\end{equation}
	A direct consequence of \eqref{eq:optimality-l-rs} is that $p^*= B \phi^*$, $b^*=q^*$ and $\mu^*=\nu^*=\eta^*$.
	
	With $\left ( \phi^*, q^*, p^*, b^*, \nu^*, \eta^* \right )$ fixed, $\mathcal{L}_{r,s}$ is concave with respect to $\mu$, and therefore $\mu^*$ verifies
	\begin{multline}
	\label{eq:optimality-mu*}
	I \left( \mu \right) - I \left( \mu^* \right) + \left \langle \mu^*-\mu, p^* -b^* \right \rangle -s \left \langle \mu^*-\mu, \mu^*-\nu^* \right \rangle \\
	-s \left \langle \mu^*-\mu, \mu^*-\eta^* \right \rangle \geq 0, \quad \forall \mu \in H,
	\end{multline}
	which implies, using $p^*= B \phi^*$ and $b^*=q^*$,
	\begin{dmath}
		F\left ( q^* \right )+ G\left ( \phi^* \right )- I \left( \mu \right)+\left \langle \mu , p^* - b^* \right \rangle +\left \langle \nu , B \phi^* - p^* \right \rangle  + \left \langle \eta , b^* - q^* \right \rangle \leq  
		F\left ( q^* \right )+ G\left ( \phi^* \right ) - I \left( \mu^* \right) +\left \langle \mu^* , p^* - b^* \right \rangle +\left \langle \nu^* , B \phi^* - p^* \right \rangle +\left \langle \eta^* , b^* - q^* \right \rangle,
		\forall \left( \mu, \nu, \eta \right) \in H^3,
	\end{dmath}
	or equivalently,
	\begin{multline}
	\label{eq:optimality-l-prime-half1}
	\mathcal{L}'\left ( \phi^*,q^*,p^*,b^*,\mu,\nu,\eta \right ) \leq \mathcal{L}'\left ( \phi^*,q^*,p^*,b^*,\mu^*,\nu^*,\eta^* \right ), \\ \forall \left( \mu, \nu, \eta \right) \in H^3.
	\end{multline}
	On the other hand, given  fixed $\left ( \mu^*, \nu^*, \eta^* , p^*, b^* \right )$, $\mathcal{L}_{r,s}$ is convex with respect to $\left ( \phi, q \right )$, and $\left ( \phi^*, q^* \right )$ satisfies
	\begin{equation}
	\label{eq:gateau}
	\left\{\begin{array}{l}
	G\left ( \phi \right )- G\left ( \phi^* \right )+\left \langle \nu^*, B \phi -B \phi^* \right  \rangle \\ \quad\quad\quad\quad + r \left \langle B \phi^* -p^*, B \phi -B \phi^* \right \rangle  \geq 0, \forall \phi \in V, \\
	F\left ( q \right )- F\left ( q^* \right )-\left \langle \eta^*, q - q^* \right \rangle \\ \quad\quad\quad\quad\quad\quad\quad + r \left \langle q^* - b^* , q - q^* \right \rangle \geq 0, \forall q \in H.
	\end{array} \right.
	\end{equation}
	which implies
	\begin{multline}
	F\left ( q^* \right )+ G\left ( \phi^* \right )+\left \langle \mu^* , B \phi^* - q^* \right \rangle  \leq F\left ( q \right )+ G\left ( \phi \right )\\ +\left \langle \mu^* , B \phi - q \right \rangle, \forall \left( \phi, q \right) \in V \times H.
	\end{multline}
	Using again $\mu^*=\nu^*=\eta^*$, we deduce that
	\begin{dmath}
		F\left ( q^* \right )+ G\left ( \phi^* \right )- I \left( \mu^* \right)+\left \langle \mu^* , p^* - b^* \right \rangle +\left \langle \nu^* , B \phi^* - p^* \right \rangle + \left \langle \eta^* , b^* - q^* \right \rangle \leq  
		F\left ( q \right )+ G\left ( \phi \right )- I \left( \mu^* \right) +\left \langle \mu^* , p - b \right \rangle +\left \langle \nu^* , B \phi - p \right \rangle +\left \langle \eta^* , b - q \right \rangle,
		\forall \left( \phi, q, p, b \right) \in V \times H^3,
	\end{dmath}
	or equivalently,
	\begin{dmath}
		\label{eq:optimality-l-prime-half2}
		\mathcal{L}'\left ( \phi^*,q^*,p^*,b^*,\mu^*,\nu^*,\eta^* \right ) \leq \mathcal{L}'\left ( \phi,q,p,b,\mu^*,\nu^*,\eta^* \right ), \forall \left( \phi, q, p, b \right) \in V \times H^3.
	\end{dmath}
	Finally, \eqref{eq:optimality-l-prime-half1} and \eqref{eq:optimality-l-prime-half2} imply that $\left ( \phi^*,q^*, p^*, b^*, \mu^*, \nu^* , \eta^* \right )$ is a saddle point of $\mathcal{L}' $. This completes the proof.
\end{proof}

\subsection{Proof Theorem 2}
\begin{proof}
	In the proof, we will use the following notation:  $\bar{\phi}^n= \phi^n -\phi^*$, $\bar{p}^n= p^n-p^*$,  $\bar{\nu}^n= \nu^n -\nu^*$, $\bar{q}^n= q^n -q^*$, $\bar{\eta}^n= \eta^n -\eta^*$, $\bar{b}^n= b^n -b^*$, and $\bar{\mu}^n= \mu^n -\mu^*$.
	We may also use the following relations without mentioning $B \phi^* = p^*$, $b^*=q^*$ and $\mu^*= \nu^* =\eta^*$.
	
	Under the assumptions of the theorem, both $\mathcal{L}_{p^n, \mu^n}$ and  $\mathcal{L}_{b^n, \mu^n}$  admit at least one  saddle point according to Theorem \ref{thm:existence-saddle-point}.
	On the other hand, computing $\mu^n$ corresponds to the minimization of a coercive, proper, convex, lower semi-continuous functional over a Hilbert space.
	By the generalized Weierstrass theorem (Theorem \ref{thm:weiestrass}), this problem admits a solution.
	This shows that the algorithm is well defined.
	
	Given that $\left ( \phi^* , q^* , p^* , b^*, \mu^* , \nu^*, \eta^*  \right )$ is a saddle point of $\mathcal{L}_{r,s}$, we have (see \eqref{eq:gateau} and use the fact that $\mu^*=\nu^*=\eta^*$)
	\begin{dmath}
		\label{eq:alg1-extremality}
		F\left ( q^n \right )+ G\left ( \phi^n \right ) +\left \langle \mu^*, p^n - b^n \right \rangle + \left \langle \nu^*, B \phi^n-p^n \right \rangle + \left \langle \eta^*, b^n -q^n \right \rangle  \geq F\left ( q^* \right )+ G\left ( \phi^* \right ) +\left \langle \mu^*, p^*-  b^* \right \rangle
	\end{dmath}
	On the other hand, by construction, $\phi^n$ and $q^n$ verify
	\begin{equation}
	\label{eq:alg1-phi-q}
	\left\{\begin{array}{l}
	G\left ( \phi \right )- G\left ( \phi^n\right )  \\ +\left \langle \nu^n + r \left( B \phi^n -p^n \right), B \phi -B \phi^n \right \rangle \geq 0,   \forall \phi \in V, \\ 
	F\left ( q \right )- F\left ( q^n\right )\\ \quad\quad\quad +\left \langle \eta^n + r \left( b^n -q^n \right), q^n -q \right \rangle \geq 0,  \forall q \in H.
	\end{array}\right.
	\end{equation}
	Taking, $\phi= \phi^*$  and $q=q^*$ in \eqref{eq:alg1-phi-q} yields
	\begin{equation}
	\label{eq:alg1-phi-q-star}
	\left\{\begin{array}{lr}
	G\left ( \phi^* \right )- G\left ( \phi^n\right ) - \left \langle \nu^n+ r \left( B \bar{\phi}^n - \bar{p}^n \right), B \bar{\phi}^n \right \rangle & \geq 0,\\ 
	F\left ( q^* \right )- F\left ( q^n\right )+\left \langle \eta^n + r \left( \bar{b}^n -\bar{q}^n \right), \bar{q}^n \right \rangle & \geq 0.
	\end{array}\right.
	\end{equation}
	By construction, $\nu^n$ and $\eta^n$ verify
	\begin{equation}
	\label{eq:alg1-nu-eta}
	\left\{\begin{array}{l}
	\nu^n= \dfrac{1}{s} \left ( B \phi^n -p^n \right )  + \mu^n,\\ 
	\eta^n= \dfrac{1}{s} \left ( b^n -q^n \right ) + \mu^n,
	\end{array}\right.
	\end{equation}
	which implies
	\begin{equation}
	\label{eq:alg1-nu-eta-prod}
	\left\{\begin{array}{lr}
	\left \langle \bar{\nu}^n, B \bar{\phi}^n - \bar{p}^n \right \rangle + s \left \langle \bar{\nu}^n,  \bar{\mu}^n - \bar{\nu}^n \right \rangle & =0,\\ 
	\left \langle \bar{\eta}^n, \bar{b}^n - \bar{q}^n \right \rangle + s \left \langle \bar{\eta}^n,  \bar{\mu}^n - \bar{\eta}^n \right \rangle & =0.
	\end{array}\right.
	\end{equation}
	Adding \eqref{eq:alg1-nu-eta-prod}  to \eqref{eq:alg1-phi-q-star} line by line gives
	\begin{multline}
	\label{eq:alg1-phi-nu-1}
	G\left ( \phi^* \right )- G\left ( \phi^n\right ) - \left \langle \nu^*, B\phi^n - p^n \right \rangle -\left \langle \nu^n , \bar{p}^n \right \rangle \\- r \left \langle B \bar{\phi}^n - \bar{p}^n , B \bar{\phi}^n \right \rangle 
	+ s \left \langle \bar{\nu}^n,  \bar{\mu}^n - \bar{\nu}^n \right \rangle\geq 0, 
	\end{multline}
	\begin{multline}
	\label{eq:alg1-phi-nu-2}
	F\left ( q^* \right )- F\left ( q^n\right ) -\left \langle \eta^*, b^n - q^n \right \rangle + \left \langle \eta^n, \bar{b}^n \right \rangle \\ + r \left \langle \bar{b}^n -\bar{q}^n , \bar{q}^n \right \rangle +  s \left \langle \bar{\eta}^n,  \bar{\mu}^n - \bar{\eta}^n \right \rangle\geq 0.
	\end{multline}
	Summing the two inequalities \eqref{eq:alg1-phi-nu-1}, \eqref{eq:alg1-phi-nu-2} along with \eqref{eq:alg1-extremality} results in
	\begin{dmath}
		-\left \langle \nu^n , \bar{p}^n \right \rangle- r \left \langle B \bar{\phi}^n - \bar{p}^n , B \bar{\phi}^n \right \rangle + s \left \langle \bar{\nu}^n,  \bar{\mu}^n - \bar{\nu}^n \right \rangle +\left \langle \eta^n, \bar{b}^n \right \rangle 
		+ r \left \langle \bar{b}^n -\bar{q}^n , \bar{q}^n \right \rangle   +  s \left \langle \bar{\eta}^n,  \bar{\mu}^n - \bar{\eta}^n \right \rangle + \left \langle \mu^*, \bar{p}^n-\bar{b}^n \right \rangle\geq 0.
	\end{dmath}
	Using again \eqref{eq:alg1-nu-eta}, we get
	\begin{dmath}
		\label{eq:alg1-sum}
		-\left \langle B \bar{\phi}^n - \bar{p}^n , r B \bar{\phi}^n + \frac{1}{s} \bar{p}^n \right \rangle +  \left \langle \bar{b}^n -\bar{q}^n , r \bar{q}^n + \frac{1}{s} \bar{b}^n \right \rangle + s \left \langle \bar{\nu}^n,  \bar{\mu}^n - \bar{\nu}^n \right \rangle + s \left \langle \bar{\eta}^n,  \bar{\mu}^n - \bar{\eta}^n \right \rangle - \left \langle \bar{\mu}^n, \bar{p}^n-\bar{b}^n \right \rangle\geq 0.
	\end{dmath}
	On the other hand, $\mu^n$ verifies by construction
	\begin{dmath}
		\label{eq:alg1-mu}
		I\left ( \mu \right) - I \left( \mu^n \right) + \left \langle \mu^n -\mu, p^n - b^n \right \rangle - s \left \langle \mu^n - \mu , \mu^n - \nu^{n-1} \right \rangle 
		- s \left \langle \mu^n - \mu , \mu^n - \eta^{n-1} \right \rangle \geq 0, \quad \forall \mu \in H.
	\end{dmath}
	Taking $\mu= \mu^*$, we obtain
	\begin{dmath}
		\label{eq:alg1-mu-star}
		I\left ( \mu^* \right) - I \left( \mu^n \right) + \left \langle \bar{\mu}^n , p^n - b^n \right \rangle - s \left \langle \bar{\mu}^n , \bar{\mu}^n - \bar{\nu}^{n-1} \right \rangle - s \left \langle \bar{\mu}^n , \bar{\mu}^n - \bar{\eta}^{n-1} \right \rangle \geq 0.
	\end{dmath}
	Summing \eqref{eq:alg1-mu-star} and \eqref{eq:alg1-sum} and rearranging yields
	\begin{dmath}
		\label{eq:alg1-before-cauchy}
		I\left ( \mu^* \right) - I \left( \mu^n \right) -\left \langle B \bar{\phi}^n - \bar{p}^n , r B \bar{\phi}^n + \frac{1}{s} \bar{p}^n \right \rangle \\ +  \left \langle \bar{b}^n -\bar{q}^n , r \bar{q}^n + \frac{1}{s} \bar{b}^n \right \rangle 
		- s \left | \mu^n -\nu^n \right |^2  - s \left | \mu^n -\eta^n \right |^2  + s \left \langle \bar{\mu}^n,  \bar{\nu}^{n-1} - \bar{\nu}^n \right \rangle 
		+ s \left \langle \bar{\mu}^n,  \bar{\eta}^{n-1} - \bar{\eta}^n \right \rangle + \left \langle \bar{\mu}^n, p^*-b^* \right \rangle\geq 0.
	\end{dmath}
	From \eqref{eq:alg1-nu-eta}, we deduce 
	\begin{equation}
	\bar{\mu}^n= \bar{\nu}^n - \frac{1}{s}\left ( B \bar{\phi}^n - \bar{p}^n \right ) = \bar{\eta}^n - \frac{1}{s}\left ( \bar{b}^n - \bar{q}^n \right ),
	\end{equation}
	which gives
	\begin{multline}
	s \left \langle \bar{\mu}^n,  \bar{\nu}^{n-1} - \bar{\nu}^n \right \rangle =  s \left \langle \bar{\nu}^n,  \bar{\nu}^{n-1} - \bar{\nu}^n \right \rangle \\ -   \left \langle B\bar{\phi}^n - \bar{p}^n,  \bar{\nu}^{n-1} - \bar{\nu}^n \right \rangle
	\end{multline}
	\begin{multline}
	s \left \langle \bar{\mu}^n,  \bar{\eta}^{n-1} - \bar{\eta}^n \right \rangle =  s \left \langle \bar{\eta}^n,  \bar{\eta}^{n-1} - \bar{\eta}^n \right \rangle \\ -  \left \langle \bar{b}^n - \bar{q}^n,  \bar{\eta}^{n-1} - \bar{\eta}^n \right \rangle.
	\end{multline}
	Using the identity
	\begin{dmath}
		\left \langle \bar{\nu}^n,  \bar{\nu}^{n-1} - \bar{\nu}^n \right \rangle =  \frac{1}{2} \left ( \left | \bar{\nu}^{n-1} \right |^2 - \left | \bar{\nu}^n \right |^2 - \left | \bar{\nu}^n -\bar{\nu}^{n-1} \right |^2 \right )
	\end{dmath}
	and its analogous for $\left \langle \bar{\eta}^n,  \bar{\eta}^{n-1} - \bar{\eta}^n \right \rangle$, and from the Cauchy-Shwartz inequality
	\begin{multline}
	-   \left \langle B\bar{\phi}^n - \bar{p}^n,  \bar{\nu}^{n-1} - \bar{\nu}^n \right \rangle \leq \frac{1}{2 \alpha}\left |  B\bar{\phi}^n - \bar{p}^n \right |^2 \\ + \frac{\alpha}{2} \left | \bar{\nu}^n - \bar{\nu}^{n-1} \right |^2, \quad \forall \alpha > 0 ,
	\end{multline}
	and
	\begin{multline}
	-   \left \langle  \bar{b}^n - \bar{q}^n,  \bar{\eta}^{n-1} - \bar{\eta}^n \right \rangle \leq \frac{1}{2 \beta}\left |   \bar{b}^n - \bar{q}^n \right |^2 \\ + \frac{\beta}{2} \left | \bar{\eta}^n - \bar{\eta}^{n-1} \right |^2, \quad \forall \beta > 0 ,
	\end{multline}
	we obtain,
	\begin{multline}
	\label{eq:alg1-cauchy-nu}
	s \left \langle \bar{\mu}^n,  \bar{\nu}^{n-1} - \bar{\nu}^n \right \rangle \leq \frac{s}{2} \left ( \left | \bar{\nu}^{n-1} \right |^2 - \left | \bar{\nu}^n \right |^2  \right ) \\ +\frac{1}{2 \alpha} \left |  B\bar{\phi}^n - \bar{p}^n \right |^2  + \left ( \frac{\alpha}{2} - \frac{s}{2} \right ) \left | \bar{\nu}^n - \bar{\nu}^{n-1} \right |^2,
	\end{multline}
	\begin{multline}
	\label{eq:alg1-cauchy-eta}
	s \left \langle \bar{\mu}^n,  \bar{\eta}^{n-1} - \bar{\eta}^n \right \rangle \leq \frac{s}{2} \left ( \left | \bar{\eta}^{n-1} \right |^2 - \left | \bar{\eta}^n \right |^2  \right )\\ +\frac{1}{2 \beta}\left |  \bar{b}^n - \bar{q}^n \right |^2 
	+ \left ( \frac{\beta}{2} - \frac{s}{2} \right ) \left | \bar{\eta}^n - \bar{\eta}^{n-1} \right |^2.
	\end{multline}
	Replacing \eqref{eq:alg1-cauchy-nu} and \eqref{eq:alg1-cauchy-eta} in \eqref{eq:alg1-before-cauchy}, gives
	\begin{dmath}
		\label{eq:alg1-before-p-b-before-mu-nu}
		-\left \langle B \bar{\phi}^n - \bar{p}^n , r B \bar{\phi}^n + \frac{1}{s} \bar{p}^n \right \rangle +  \left \langle \bar{b}^n -\bar{q}^n , r \bar{q}^n + \frac{1}{s} \bar{b}^n \right \rangle   - s \left | \mu^n -\nu^n \right |^2  - s \left | \mu^n -\eta^n \right |^2
		+ \frac{s}{2} \left ( \left | \bar{\nu}^{n-1} \right |^2 - \left | \bar{\nu}^n \right |^2  \right )  +\frac{1}{2 \alpha} \left |  B\bar{\phi}^n - \bar{p}^n \right |^2  + \frac{\alpha - s}{2} \left | \bar{\nu}^n - \bar{\nu}^{n-1} \right |^2 
		+\frac{s}{2} \left ( \left | \bar{\eta}^{n-1} \right |^2 - \left | \bar{\eta}^n \right |^2  \right ) +\frac{1}{2 \beta}\left |  \bar{b}^n - \bar{q}^n \right |^2 \\+ \frac{\beta - s}{2} \left | \bar{\eta}^n - \bar{\eta}^{n-1} \right |^2 
		+I\left ( \mu^* \right) - I \left( \mu^n \right)+ \left \langle \bar{\mu}^n, p^*-b^* \right \rangle \geq 0.
	\end{dmath}
	From \eqref{eq:alg1-nu-eta} we can see that
	\begin{align}
	s\left | \bar{\mu}^n - \bar{\nu}^n \right |^2 &= \frac{1}{s} \left | B \bar{\phi}^n - p^n \right |^2 \\
	s\left | \bar{\mu}^n - \bar{\eta}^n \right |^2 &= \frac{1}{s} \left | b^n - q^n \right |^2
	\end{align}
	which, when used in \eqref{eq:alg1-before-p-b-before-mu-nu} results in
	\begin{dmath}
		\label{eq:alg1-before-p-b}
		-\left \langle B \bar{\phi}^n - \bar{p}^n , r B \bar{\phi}^n + \frac{1}{s} \bar{p}^n \right \rangle +  \left \langle \bar{b}^n -\bar{q}^n , r \bar{q}^n + \frac{1}{s} \bar{b}^n \right \rangle  + \frac{s}{2} \left ( \left | \bar{\nu}^{n-1} \right |^2 - \left | \bar{\nu}^n \right |^2  \right ) -\left ( \frac{1}{s} - \frac{1}{2 \alpha} \right ) \left |  B\bar{\phi}^n - \bar{p}^n \right |^2  + \left ( \frac{\alpha}{2} - \frac{s}{2} \right ) \left | \bar{\nu}^n - \bar{\nu}^{n-1} \right |^2 + \frac{s}{2} \left ( \left | \bar{\eta}^{n-1} \right |^2 - \left | \bar{\eta}^n \right |^2  \right ) 
		-\left ( \frac{1}{s} - \frac{1}{2 \beta} \right )\left |  \bar{b}^n - \bar{q}^n \right |^2 + \left ( \frac{\beta}{2} - \frac{s}{2} \right ) \left | \bar{\eta}^n - \bar{\eta}^{n-1} \right |^2 
		+I\left ( \mu^* \right) - I \left( \mu^n \right)+ \left \langle \bar{\mu}^n, p^*-b^* \right \rangle \geq 0.
	\end{dmath}
	Taking the square norm of the two sides of each the assignments in Step 2 and 4 of Algorithm 1 gives
	\begin{dmath}
		\label{eq:alg1-p-square}
		\left | \bar{p}^n \right |^2 - \left | \bar{p}^{n+1} \right |^2= - \rho^2 \left ( r+\frac{1}{s} \right )^2 \left | B \bar{\phi}^n - \bar{p}^n \right |^2 - 2 \rho \left ( r+\frac{1}{s} \right ) \left \langle \bar{p}^n, B\bar{\phi}^n - \bar{p}^n \right \rangle,
	\end{dmath}
	and
	\begin{dmath}
		\label{eq:alg1-b-square}
		\left | \bar{b}^n \right |^2 - \left | \bar{b}^{n+1} \right |^2= - \rho^2 \left ( r+\frac{1}{s} \right )^2 \left | \bar{b}^n - \bar{q}^n \right |^2  + 2 \rho \left ( r+\frac{1}{s} \right ) \left \langle \bar{b}^n, \bar{b}^n - \bar{q}^n \right \rangle. 
	\end{dmath}
	Multiplying \eqref{eq:alg1-before-p-b} by $2 \rho$ and summing  with \eqref{eq:alg1-p-square} and \eqref{eq:alg1-b-square} results in 
	\begin{dmath}
		\label{eq:alg1-decreasing}
		\left ( \left | \bar{p}^{n} \right |^2  + \left | \bar{b}^n \right |^2 + \rho s \left | \bar{\nu}^{n-1} \right |^2 + \rho s \left | \bar{\eta}^{n-1} \right |^2 \right ) 
		- \left ( \left | \bar{p}^{n+1} \right |^2 + \left | \bar{b}^{n+1} \right |^2  + \rho s \left | \bar{\nu}^n \right |^2 + \rho s \left | \bar{\eta}^n \right |^2 \right ) \geq  
		\left ( 2 \rho \left ( r + \frac{1}{s} - \frac{1}{2\alpha} \right )- \rho^2 \left ( r+\frac{1}{s} \right )^2 \right ) \left | B \bar{\phi}^n - \bar{p}^n \right |^2 + \left ( \frac{s}{2} - \frac{\alpha}{2} \right ) \left | \bar{\nu}^n - \bar{\nu}^{n-1} \right |^2 
		+ \left ( 2 \rho \left ( r + \frac{1}{s} - \frac{1}{2\beta} \right )- \rho^2 \left ( r+\frac{1}{s} \right )^2 \right ) \left | \bar{b}^n - \bar{q}^n \right |^2 + \left ( \frac{s}{2} - \frac{\beta}{2} \right ) \left | \bar{\eta}^n - \bar{\eta}^{n-1} \right |^2  
		+ 2 \rho \left ( I\left ( \mu^n \right) - I \left( \mu^* \right)  - \left \langle \bar{\mu}^n, p^*-b^* \right \rangle \right ).
	\end{dmath}
	For any $\rho$ satisfying
	\begin{equation}
	0 < \rho < \frac{2rs^2 +s}{ \left ( 1+ rs \right )^2},
	\end{equation}
	we can choose $\alpha$ and $\beta$ such that
	\begin{equation}
	0< \frac{1}{2\left ( r+\frac{1}{s} \right ) - \rho \left ( r+\frac{1}{s} \right )^2} < \alpha, \beta < s,
	\end{equation}
	which in light of \eqref{eq:alg1-decreasing} implies that the sequence 
	\begin{equation}
	\label{eq:alg1:p-b-nu-eta-decreasing}
	\left | \bar{p}^{n} \right |^2  + \left | \bar{b}^n \right |^2 + \rho s \left | \bar{\nu}^{n-1} \right |^2 + \rho s \left | \bar{\eta}^{n-1} \right |^2 \text{ is decreasing,}
	\end{equation}
	and 
	\begin{dmath}
		\left | B \bar{\phi}^n - \bar{p}^n \right |^2  + \left | \bar{b}^n - \bar{q}^n \right |^2 + \left ( \frac{s}{2} - \frac{\alpha}{2} \right ) \left | \bar{\nu}^n - \bar{\nu}^{n-1} \right |^2 + \left ( \frac{s}{2} - \frac{\beta}{2} \right ) \left | \bar{\eta}^n - \bar{\eta}^{n-1} \right |^2 
		+ 2 \rho \left ( I\left ( \mu^n \right) - I \left( \mu^* \right)  - \left \langle \bar{\mu}^n, p^*-b^* \right \rangle \right )
	\end{dmath}
	converges to zero (remember from \eqref{eq:optimality-mu*} that $I\left ( \mu^n \right) - I \left( \mu^* \right)  - \left \langle \bar{\mu}^n, p^*-b^* \right \rangle \geq 0, \forall \mu^n \in H$).
	It follows that 
	\begin{equation}
	\label{eq:alg1-bounded}
	\left | \bar{p}^{n} \right |^2,  \left | \bar{b}^n \right |^2,  \left | \bar{\nu}^n \right |^2, \left | \bar{\eta}^n \right |^2, \left | B\bar{\phi}^{n} \right |^2,  \left | \bar{q}^n \right |^2, \left | \bar{\mu}^{n} \right |^2, \text{are bounded},
	\end{equation}
	and 
	\begin{dmath}
		\label{eq:alg1-convergence-to-zero}
		\lim_{n \to \infty} \left | B \phi^n - p^n \right | = \lim_{n \to \infty} \left | b^n - q^n \right | = \lim_{n \to \infty} \left |\mu^n - \nu^n \right | = \lim_{n \to \infty} \left |\mu^n - \eta^n \right | 
		= \lim_{n \to \infty} \left |\nu^n - \nu^{n-1} \right | = \lim_{n \to \infty} \left |\eta^n - \eta^{n-1} \right |= \lim_{n \to \infty} I\left ( \mu^n \right) - I \left( \mu^* \right)  - \left \langle \bar{\mu}^n, p^*-b^* \right \rangle=0.
	\end{dmath}
	Using the assignments in Step 2 and 4 of Algorithm 1, this implies:
	\begin{equation}
	\label{eq:alg1-conv-p-conv-q}
	\lim_{n \to \infty} \left | p^{n+1} - p^n \right | = \lim_{n \to \infty} \left | q^{n+1} - q^n \right | = 0
	\end{equation}
	
	From \eqref{eq:alg1-phi-nu-1} and \eqref{eq:alg1-phi-nu-2}, we can deduce that
	\begin{dmath}
		\label{eq:alg1-phi-nu-summed}
		F\left ( q^* \right ) + G\left ( \phi^* \right )\geq F\left ( q^n\right )+ G\left ( \phi^n\right ) + \left \langle \nu^*, B\phi^n - p^n \right \rangle + \left \langle \eta^*, b^n - q^n \right \rangle 
		+ r \left \langle B \bar{\phi}^n - \bar{p}^n , B \bar{\phi}^n \right \rangle - r \left \langle \bar{b}^n -\bar{q}^n , \bar{q}^n \right \rangle - s \left \langle \bar{\nu}^n,  \bar{\mu}^n - \bar{\nu}^n \right \rangle - s \left \langle \bar{\eta}^n,  \bar{\mu}^n - \bar{\eta}^n \right \rangle 
		+ \left \langle \nu^n , \bar{p}^n \right \rangle  - \left \langle \eta^n, \bar{b}^n \right \rangle. 
	\end{dmath}
	The last two terms can be written as
	\begin{dmath}
		\left \langle \nu^n , \bar{p}^n \right \rangle  - \left \langle \eta^n, \bar{b}^n \right \rangle = \left \langle \nu^n - \mu^n, \bar{p}^n \right \rangle  - \left \langle \eta^n - \mu^n , \bar{b}^n \right \rangle + \left \langle \mu^n , \bar{p}^n - \bar{b}^n \right \rangle,
	\end{dmath}
	which, when replaced in \eqref{eq:alg1-phi-nu-summed} gives:
	\begin{dmath}
		\label{eq:alg1-limit-up}
		F\left ( q^* \right ) + G\left ( \phi^* \right ) + \left \langle \mu^n , p^* - b^* \right \rangle \geq F\left ( q^n\right )+ G\left ( \phi^n\right )  + \left \langle \mu^n , p^n - b^n \right \rangle 
		+ \left \langle \nu^*, B\phi^n - p^n \right \rangle + \left \langle \eta^*, b^n - q^n \right \rangle + r \left \langle B \bar{\phi}^n - \bar{p}^n , B \bar{\phi}^n \right \rangle - r \left \langle \bar{b}^n -\bar{q}^n , \bar{q}^n \right \rangle 
		- s \left \langle \bar{\nu}^n,  \bar{\mu}^n - \bar{\nu}^n \right \rangle - s \left \langle \bar{\eta}^n,  \bar{\mu}^n - \bar{\eta}^n \right \rangle + \left \langle \nu^n - \mu^n, \bar{p}^n \right \rangle  - \left \langle \eta^n - \mu^n , \bar{b}^n \right \rangle. 
	\end{dmath}
	On the other hand, \eqref{eq:alg1-extremality} gives us
	\begin{dmath}
		\label{eq:lagrangian-star-compared-to-n}
		F\left ( q^n \right )+ G\left ( \phi^n \right )+\left \langle \mu^*, p^n - b^n \right \rangle + \left \langle \nu^*, B \phi^n-p^n \right \rangle  
		+ \left \langle \eta^*, b^n -q^n \right \rangle  \geq F\left ( q^* \right )+ G\left ( \phi^* \right ) +\left \langle \mu^*, p^*-  b^* \right \rangle.
	\end{dmath}
	Adding \eqref{eq:alg1-mu-star} to \eqref{eq:lagrangian-star-compared-to-n} results in
	\begin{dmath}
		\label{eq:alg1-limit-low}
		F\left ( q^n \right )+ G\left ( \phi^n \right ) - I \left( \mu^n \right)+\left \langle \mu^n, p^n - b^n \right \rangle + \left \langle \nu^*, B \phi^n-p^n \right \rangle + \left \langle \eta^*, b^n -q^n \right \rangle
		- s \left \langle \bar{\mu}^n , \mu^n - \nu^n \right \rangle  - s \left \langle \bar{\mu}^n , \mu^n - \eta^n \right \rangle - s \left \langle \bar{\mu}^n , \nu^n - \nu^{n-1} \right \rangle  - s \left \langle \bar{\mu}^n , \eta^n - \eta^{n-1} \right \rangle 
		\geq F\left ( q^* \right )+ G\left ( \phi^* \right ) - I\left ( \mu^* \right)+\left \langle \mu^*, p^*-  b^* \right \rangle.
	\end{dmath}
	From  \eqref{eq:alg1-limit-up}, \eqref{eq:alg1-bounded} and \eqref{eq:alg1-convergence-to-zero}, it follows that
	\begin{dmath}
		\label{eq:alg1-lim-sup}
		\limsup \left [ F\left ( q^* \right ) + G\left ( \phi^* \right ) - I \left( \mu^n \right )+ \left \langle \mu^n , p^* - b^* \right \rangle \right ] 
		\geq \limsup\left [  F\left ( q^n\right )+ G\left ( \phi^n\right ) - I \left( \mu^n \right ) + \left \langle \mu^n , p^n - b^n \right \rangle \right ].
	\end{dmath}
	From  \eqref{eq:alg1-limit-low}, \eqref{eq:alg1-bounded} and \eqref{eq:alg1-convergence-to-zero}, it follows that
	\begin{dmath}
		\label{eq:alg1-lim-inf}
		\liminf \left [ F\left ( q^n \right )+ G\left ( \phi^n \right ) - I \left( \mu^n \right )+\left \langle \mu^n, p^n - b^n \right \rangle \right ]
		\geq F\left ( q^* \right )+ G\left ( \phi^* \right )  - I \left( \mu^* \right )+\left \langle \mu^*, p^*-  b^* \right \rangle.
	\end{dmath}
	From the fact that $I\left ( \mu^n \right) - I \left( \mu^* \right)  - \left \langle \bar{\mu}^n, p^*-b^* \right \rangle \geq 0, \forall \mu^n \in H$, we have
	\begin{dmath}
		F\left ( q^* \right )+ G\left ( \phi^* \right )  - I \left( \mu^* \right )+\left \langle \mu^*, p^*-  b^* \right \rangle 
		\geq \limsup \left [ F\left ( q^* \right ) + G\left ( \phi^* \right ) - I \left( \mu^n \right )+ \left \langle \mu^n , p^* - b^* \right \rangle \right ],
	\end{dmath}
	which when combined with \eqref{eq:alg1-lim-sup} and \eqref{eq:alg1-lim-inf} implies that
	\begin{multline}
	\label{eq:alg1-convergence-lagrangian-value}
	\lim_{n \to \infty} \left [ F\left ( q^n \right ) + G\left ( \phi^n \right )- I \left( \mu^n \right)+\left \langle \mu^n, p^n - b^n \right \rangle \right ] \\
	= F\left ( q^* \right )+ G\left ( \phi^* \right ) - I \left( \mu^* \right) +\left \langle \mu^*, p^*-  b^* \right \rangle.
	\end{multline}
	Assume that $G$ is uniformly convex over the bounded sets of $V$ and let $\nu^{n+1/2}$ be such that $-B^t\nu^{n+1/2} \in \partial G\left ( \left( \phi^n+\phi^* \right ) /2 \right ) $.
	Then
	\begin{align*}
	& F\left ( q^n \right ) - F\left ( q^*\right )+\left \langle \eta^*,q^*-q^n \right \rangle \geq 0, \\
	& G\left ( \phi^n \right ) - G\left ( \frac{\phi^n+\phi^*}{2} \right )+\left \langle \nu^{n+1/2}, B\phi^n - \frac{B\phi^n+B\phi^*}{2} \right \rangle \geq 0 ,
	\end{align*}
	which imply
	\begin{dmath}
		F\left ( q^n \right ) + G\left ( \phi^n \right ) + \left \langle \mu^*,B \phi^n-q^n \right \rangle \geq  G\left ( \frac{\phi^n+\phi^*}{2} \right ) - G\left ( \phi^* \right ) 
		+ \left \langle \nu^*, \frac{B\phi^n +B\phi^*}{2}-B\phi^* \right \rangle + \left \langle \nu^* - \nu^{n+1/2},\frac{B\phi^n +B\phi^*}{2}-B\phi^* \right \rangle
		+F\left ( q^*\right ) + G\left ( \phi^* \right )+\left \langle \mu^*, B\phi^*-q^* \right \rangle,
	\end{dmath}
	Since
	\begin{dmath}
		G\left ( \frac{\phi^n+\phi^*}{2} \right ) - G\left ( \phi^* \right ) + \left \langle \nu^*, \frac{B\phi^n +B\phi^*}{2}-B\phi^* \right \rangle \geq 0,
	\end{dmath}
	it follows that
	\begin{dmath}
		F\left ( q^n \right ) + G\left ( \phi^n \right ) + \left \langle \mu^*,B \phi^n-q^n \right \rangle \geq  \left \langle \nu^* - \nu^{n+1/2},\frac{B\phi^n +B\phi^*}{2}-B\phi^* \right \rangle
		+F\left ( q^*\right ) + G\left ( \phi^* \right )+\left \langle \mu^*, B\phi^*-q^* \right \rangle.
	\end{dmath}
	From Theorem 1 (more precisely \eqref{eq:gateau}), we have:
	\begin{equation}
	-B^t \nu^* \in \partial G \left ( \phi^*\right ) \text{ and } \eta^* \in \partial F \left ( q^*\right )
	\end{equation}
	By the definition of uniform convexity, we obtain
	\begin{dmath}
		F\left ( q^n \right ) + G\left ( \phi^n \right ) + \left \langle \mu^*,B \phi^n-q^n \right \rangle \geq  \Psi \left(  \frac{\left | \phi^n-\phi^*\right | }{2}  \right )
		+F\left ( q^*\right ) + G\left ( \phi^* \right )+\left \langle \mu^*, B\phi^*-q^* \right \rangle.
	\end{dmath}
	By \eqref{eq:alg1-mu-star},
	\begin{dmath}
		F\left(q^n \right) + G\left(\phi^n \right) - I\left(\mu^n \right) + \left\langle \mu^n, p^n-b^n \right \rangle + \left\langle \mu^*, B \phi^n-p^n \right \rangle
		+ \left \langle \mu^* , b^n - q^n \right \rangle - s \left \langle \bar{\mu}^n , \bar{\mu}^n - \bar{\nu}^{n-1} \right \rangle - s \left \langle \bar{\mu}^n , \bar{\mu}^n - \bar{\eta}^{n-1} \right \rangle
		\geq  \Psi \left(  \frac{\left | \phi^n-\phi^*\right | }{2}  \right )+F\left ( q^*\right ) + G\left ( \phi^* \right )-I\left ( \mu^* \right)+\left \langle \mu^*, B\phi^*-q^* \right \rangle.
	\end{dmath}
	From \eqref{eq:alg1-bounded}, \eqref{eq:alg1-convergence-to-zero} and \eqref{eq:alg1-convergence-lagrangian-value}, it follows that 
	\begin{equation}
	\lim_{n \to \infty} \Psi \left(  \frac{\left | \phi^n-\phi^*\right | }{2}  \right ) =0,
	\end{equation}
	which, by the definition of uniform convexity, implies $\phi^n \to \phi^*$.
	From \eqref{eq:alg1-convergence-to-zero} and \eqref{eq:alg1-conv-p-conv-q}, we also have that $B \phi^n$ converges.
	Since $B$ is closed range and injective, it follows that $B \phi^n \to B \phi^*$ and hence $p^n \to B \phi^*$ .
	A similar argument for $q$ yields $q^n \to q^*$. It follows then from \eqref{eq:alg1-convergence-to-zero} that $b^n \to q^*$.
	
	Now, we apply Opial's lemma  (Lemma \ref{lem:opial}) to prove the last statement of the proposition.
	Let $W=\left \{ \hat{\mu}\in H,  \left ( \phi^*, q^*, p^*, b^*, \hat{\mu}, \hat{\mu}, \hat{\mu} \right )  \text{ is a saddle point of } \mathcal{L}_{r,s}\right \}$.
	$W$ is clearly closed.
	Since the sequence $\left ( \mu^n \right)_n$ is bounded, we can extract  a subsequence  $\left (\mu^{n_k}\right )_{n_k}$ that converges weakly to a point $\hat{\mu}$.
	
	Then by \eqref{eq:alg1-convergence-to-zero}, passing to the limit superior in \eqref{eq:alg1-mu} and the fact that $I$ is weakly lower semi-continuous (since $I$ is convex and lower semi-continuous), we have that:
	\begin{equation}
	- I \left( \hat{\mu} \right) + \left \langle \hat{\mu}, p^* - b^* \right \rangle \geq -I\left ( \mu \right) +  \left \langle \mu, p^* - b^* \right \rangle,  \forall \mu \in H.
	\end{equation}
	
	Also, by \eqref{eq:alg1-convergence-to-zero}, \eqref{eq:alg1-convergence-lagrangian-value}, summing and passing to the limit in \eqref{eq:alg1-phi-q}, we have that:
	\begin{multline}
	G\left ( \phi \right )  +  F\left ( q \right ) + \left \langle \hat{\mu}, B \phi - q \right \rangle  \geq  
	G\left ( \phi^*\right )+ F\left ( q^*\right ) \\ +\left \langle \hat{\mu}, B \phi^* - q^* \right \rangle, \forall  \phi \in V, q \in H.
	\end{multline}
	The two last results imply that $\hat{\mu} \in W$.
	
	From \eqref{eq:alg1:p-b-nu-eta-decreasing} and what we have just proved for $p^n$ and $b^n$, we deduce that $\left | \nu^n -\mu^* \right |^2+\left | \eta^n -\mu^* \right |^2$ converges.
	It follows then from \eqref{eq:alg1-bounded} and \eqref{eq:alg1-convergence-to-zero} that $\left | \mu^n - \mu^* \right |$ converges.
	By Opial's lemma, $\mu^n$ converges weakly to $\hat{\mu} \in W$, and by \eqref{eq:alg1-convergence-to-zero} so do $\nu^n$ and $\eta^n$.
	
	Finally, from Theorem 1 (more precisely \eqref{eq:optimality-mu*}), we have:
	\begin{equation}
	B \phi^* - q^*\in \partial I \left ( \mu^*\right ).
	\end{equation}
	Therefore, if $I$ is uniformly convex on the subsets of $H$, we can use a similar argument to the one used for $\phi^n$ and $q^n$ to prove the strong convergence $\mu^n$ as well as $\nu^n$ and $\eta^n$ to $\mu^*$.
	This completes the proof.
\end{proof}

\subsection{Proof of Lemma 1}
\begin{proof}
	The existence of a saddle point of $\mathcal{L}_{p^n, \mu^n}$ has already been shown in the proof of convergence of Algorithm 1.
	Since $\left ( \phi^n, \nu^n \right ) $ is a saddle point of $\mathcal{L}_{p^n, \mu^n}$, we have 
	\begin{dmath}
		\label{eq:alg-phi-nu-saddle}
		G\left ( \phi \right ) - G\left ( \phi^n \right ) + \left \langle \nu^n, B \phi - B \phi^n  \right \rangle + r \left \langle B \phi^n - p^n, B \phi - B \phi^n \right \rangle \geq 0, \quad \forall \phi \in V.
	\end{dmath}
	On the other hand, $\phi^n_{k}$ verifies by construction
	\begin{dmath}
		\label{eq:alg-phi-nu-phi-n-k}
		G\left ( \phi \right ) - G\left ( \phi_k^n \right ) + \left \langle \nu_{k-1}^n, B \phi - B \phi_k^n  \right \rangle + r \left \langle B \phi_k^n - p^n, B \phi - B \phi_k^n \right \rangle \geq 0, \quad \forall \phi \in V.
	\end{dmath}
	Taking $\phi= \phi_k^n$ in \eqref{eq:alg-phi-nu-saddle} and $\phi= \phi^n$ in \eqref{eq:alg-phi-nu-phi-n-k}, and suming the two inequalities yields
	\begin{equation}
	\left \langle \nu_{k-1}^n - \nu^n,  B \phi^n_k - B \phi^n \right \rangle + r \left | B \phi_{k}^n - B \phi^n \right |^2 \ \leq 0.
	\end{equation}
	Putting $\bar{\nu}_{k-1}^n = \nu_{k-1}^n - \nu^n$ and $\bar{\phi}_k^n= \phi_k^n - \phi^n$, we then have
	\begin{equation}
	\label{eq:alg-phi-nu-phi-bar-nu-bar}
	\left \langle \bar{\nu}_{k-1}^n ,  B \bar{\phi}^n_k \right \rangle + r \left | B \bar{\phi}^n_{k} \right |^2 \ \leq 0.
	\end{equation}
	On the other hand, by assumption $\nu^n$ verifies
	\begin{equation}
	\mathcal{L}_{p^n,\mu^n} \left ( \phi^n, \nu \right ) \leq \mathcal{L}_{p^n,\mu^n} \left ( \phi^n, \nu^n \right ) ,
	\quad \forall  \nu \in H, 
	\end{equation}
	which implies that
	\begin{equation}
	\label{eq:alg-phi-nu-nu-saddle}
	\nu^n= \frac{B \phi^n -p^n}{s} + \mu^n.
	\end{equation}
	The update in Step 2 of Algorithm 2 implies
	\begin{dmath}
		\bar{\nu}^n_k= \bar{\nu}_{k-1}^n+ \rho_{\nu} \left ( B \bar{\phi}_{k}^n + B \phi^n- p^n - s \left ( \bar{\nu}_{k-1}^n +\nu^n - \mu^n \right ) \right ),
	\end{dmath}
	which together with \eqref{eq:alg-phi-nu-nu-saddle} results in
	\begin{equation}
	\bar{\nu}^n_k= \left ( 1- s \rho_{\nu} \right )\bar{\nu}_{k-1}^n+ \rho_{\nu} B \bar{\phi}_{k}^n .
	\end{equation}
	Taking the square norm of the two sides, we obtain 
	\begin{dmath}
		\left | \bar{\nu}^n_k \right |^2  = \left ( 1- s \rho_{\nu} \right )^2 \left | \bar{\nu}_{k-1}^n \right |^2 + \rho_{\nu}^2 \left | B \bar{\phi}_{k}^n \right |^2 + 2 \rho_{\nu} \left ( 1- s \rho_{\nu} \right ) \left \langle \bar{\nu}_{k-1}^n, B \bar{\phi}_{k}^n \right \rangle.
	\end{dmath}
	From \eqref{eq:alg-phi-nu-phi-bar-nu-bar}, 
	\begin{dmath}
		\left | \bar{\nu}^n_{k-1} \right |^2-\left | \bar{\nu}^n_{k} \right |^2 \geq \left ( 2s\rho_{\nu}-s^2\rho_{\nu}^2 \right )\left | \bar{\nu}^n_{k-1} \right |^2+\left ( 2r \rho_{\nu}\left ( 1-s\rho_{\nu} \right ) -\rho_{\nu}^2\right )\left | B\bar{\phi^n_k} \right |^2.
	\end{dmath}
	It follows that for
	\begin{equation}
	0 < \rho_{\nu} <\frac{2r}{2rs+1} < \frac{2}{s}, 
	\end{equation}
	$\left | \bar{\nu}^n_{k} \right |^2 \to 0$ and $\left | B\bar{\phi^n_k} \right |^2 \to 0$.
	Since $B$ has a closed range, $\phi^n \to \phi^*$, which completes the proof.
\end{proof}

\subsection{Proof of Theorem 3}
\begin{proof}
	In the proof, we will use the following notation:  $\bar{\phi}^n= \phi^n -\phi^*$, $\bar{p}^n= p^n-p^*$,  $\bar{\nu}^n= \nu^n -\nu^*$, $\bar{q}^n= q^n -q^*$, $\bar{\eta}^n= \eta^n -\eta^*$, $\bar{b}^n= b^n -b^*$, and $\bar{\mu}^n= \mu^n -\mu^*$.
	We may also use the following relations without mentioning $B \phi^* = p^*$, $b^*=q^*$ and $\mu^*= \nu^* =\eta^*$.
	
	First, Step 1, 4 and 7 correspond to the minimization of coercive, proper, convex, lower semi-continuous functionals over a Hilbert space.
	By the generalized Weierstrass theorem (Theorem \ref{thm:weiestrass}), these problems admit solutions.
	This shows that the algorithm is well defined.
	
	Given that $\left ( \phi^* , q^* , p^* , b^*, \mu^* , \nu^*, \eta^*  \right )$ is a saddle point of $\mathcal{L}_{r,s}$, we have (see \eqref{eq:gateau} and use the fact that $\mu^*=\nu^*=\eta^*$)
	\begin{dmath}
		\label{eq:alg2-extremality}
		F\left ( q^n \right )+ G\left ( \phi^n \right )+\left \langle \mu^*, p^{n-1} - b^{n-1} \right \rangle + \left \langle \nu^*, B \phi^n-p^{n-1} \right \rangle + \left \langle \eta^*, b^{n-1} -q^n \right \rangle  \geq F\left ( q^* \right )+ G\left ( \phi^* \right ) +\left \langle \mu^*, p^*-  b^* \right \rangle
	\end{dmath}
	On the other hand, by construction, $\phi^n$, $q^n$ and $\mu^n$ verify
	\begin{dmath}
		\label{eq:alg2-phi-q}
		\left\{\begin{array}{l}
			G\left ( \phi \right )- G\left ( \phi^n\right ) + \left \langle \nu^{n-1}, B \phi -B \phi^n  \right \rangle \\ \quad\quad+ r \left \langle B \phi^n -p^{n-1}, B \phi -B \phi^n \right \rangle \geq 0,  \forall \phi \in V, \\ 
			F\left ( q \right )- F\left ( q^n\right )  +\left \langle \eta^{n-1}, q^n -q \right \rangle \\ \quad\quad\quad\quad\quad + r \left \langle b^{n-1} -q^n , q^n -q \right \rangle \geq 0,   \forall q \in H, \\
			I\left ( \mu \right) - I \left( \mu^n \right)  +\left \langle \mu^n -\mu, p^{n-1} - b^{n-1}\right \rangle \\ \quad +\left \langle \mu^n -\mu,- s \left( \mu^n - \nu^{n-1} \right )\right \rangle \\ \quad\quad + \left \langle \mu^n -\mu,- s \left ( \mu^n - \eta^{n-1} \right ) \right \rangle \geq 0, \forall \mu \in H.
		\end{array}\right.
	\end{dmath}
	Taking, $\phi= \phi^*$,  $q=q^*$ and $\mu=\mu^*$ in \eqref{eq:alg2-phi-q} yields
	\begin{equation}
	\label{eq:alg2-phi-q-star}
	\left\{\begin{array}{l}
	G\left ( \phi^* \right )- G\left ( \phi^n\right ) \\ \quad\quad\quad - \left \langle \nu^{n-1} + r \left( B \bar{\phi}^n - \bar{p}^{n-1} \right), B \bar{\phi}^n \right \rangle  \geq 0,\\ 
	F\left ( q^* \right )- F\left ( q^n\right ) \\ \quad\quad\quad\quad\quad +\left \langle \eta^{n-1} + r \left( \bar{b}^{n-1} -\bar{q}^n \right), \bar{q}^n \right \rangle \geq 0, \\
	I\left ( \mu^* \right) - I \left( \mu^n \right) +\left \langle \bar{\mu}, p^{n-1} - b^{n-1}\right \rangle \\ \quad\quad\quad +\left \langle \bar{\mu}^n,- s \left( \mu^n - \nu^{n-1} \right )\right \rangle \\ \quad\quad\quad\quad\quad\quad\quad\quad + \left \langle \bar{\mu}^n,- s \left ( \mu^n - \eta^{n-1} \right ) \right \rangle \geq 0.
	\end{array}\right.
	\end{equation}
	Summing all inequalities in \eqref{eq:alg2-phi-q-star} along with \eqref{eq:alg2-extremality} yields
	\begin{dmath}
		\label{eq:alg2-sum-all-extremal-ineq}
		\left \langle \bar{\nu}^{n-1},  B \bar{\phi}^n \right \rangle + r \left \langle B \bar{\phi}^n - \bar{p}^{n-1} , B \bar{\phi}^n \right \rangle - \left \langle \bar{\eta}^{n-1},\bar{q}^n \right \rangle  - r \left \langle \bar{b}^{n-1} -\bar{q}^n , \bar{q}^n \right \rangle 
		+I\left ( \mu^n \right) - I \left( \mu^* \right) - \left \langle \bar{\mu}^n, p^{n-1} - b^{n-1} - s \left( \mu^n - \nu^{n-1} \right ) - s \left ( \mu^n - \eta^{n-1} \right ) \right \rangle  \leq 0.
	\end{dmath}
	Taking the square norm of each of the side of the assignments in Step 2, 3, 5 and 6 of Algorithm 4 respectively gives us
	\begin{multline}
	\label{eq:alg2-square-norm-p}
	\left | \bar{p}^{n-1} \right |^2-\left | \bar{p}^n \right |^2 = -\rho_r^2\left | \mu^n - \nu^{n-1} + r \left ( p^{n-1} - B\phi^n \right) \right |^2 \\
	+ 2 \rho_r \left \langle \bar{p}^{n-1}, \mu^n - \nu^{n-1} + r \left ( p^{n-1} - B\phi^n \right) \right \rangle,
	\end{multline}
	\begin{multline}
	\label{eq:alg2-square-norm-b}
	\left | \bar{b}^{n-1} \right |^2-\left | \bar{b}^n \right |^2 = -\rho_r^2\left | \eta^{n-1} - \mu^n + r \left ( b^{n-1} - q^n \right) \right |^2 \\
	+ 2 \rho_r \left \langle \bar{b}^{n-1}, \eta^{n-1} - \mu^n + r \left ( b^{n-1} - q^n \right) \right \rangle,
	\end{multline}
	\begin{multline}
	\label{eq:alg2-square-norm-nu}
	\left | \bar{\nu}^{n-1} \right |^2-\left | \bar{\nu}^n \right |^2 = -\rho_s^2\left | B\phi^n - p^{n-1} -s \left ( \nu^{n-1} -  \mu^n \right) \right |^2 \\
	- 2 \rho_s \left \langle \bar{\nu}^{n-1}, B\phi^n - p^{n-1} -s \left ( \nu^{n-1} -  \mu^n \right) \right \rangle,
	\end{multline}
	\begin{multline}
	\label{eq:alg2-square-norm-eta}
	\left | \bar{\eta}^{n-1} \right |^2-\left | \bar{\eta}^n \right |^2 = -\rho_s^2\left | b^{n-1} - q^n -s \left ( \eta^{n-1} -  \mu^n \right) \right |^2 \\
	- 2 \rho_s \left \langle \bar{\eta}^{n-1}, b^{n-1} - q^n -s \left ( \eta^{n-1} -  \mu^n \right) \right \rangle.
	\end{multline}
	Summing \eqref{eq:alg2-square-norm-p} and \eqref{eq:alg2-square-norm-nu} yields
	\begin{dmath}
		\label{eq:alg2-sum-square-norm-p-nu}
		\rho_s \left ( \left | \bar{p}^{n-1} \right |^2-\left | \bar{p}^n \right |^2 \right ) + \rho_r \left ( \left | \bar{\nu}^{n-1} \right |^2-\left | \bar{\nu}^n \right |^2 \right )= - \rho_r \rho_s \left ( \rho_r + \rho_s s^2 \right ) \left | \mu^n - \nu^{n-1} \right|^2 
		- \rho_r \rho_s \left ( \rho_r r^2 + \rho_s \right ) \left | B \phi^n - p^{n-1}  \right |^2 + 2 \rho_s \rho_r \left \langle \bar{p}^{n-1}, \mu^n - \nu^{n-1} + r \left ( p^{n-1} - B\phi^n \right) \right \rangle 
		- 2 \rho_r \rho_s \left \langle \bar{\nu}^{n-1}, B\phi^n - p^{n-1} -s \left ( \nu^{n-1} -  \mu^n \right) \right \rangle\\
		+ 2 \rho_r \rho_s \left ( \rho_r r - \rho_s s  \right ) \left \langle \mu^n - \nu^{n-1}, B\phi^n - p^{n-1} \right \rangle,
	\end{dmath}
	whereas summing \eqref{eq:alg2-square-norm-b} and \eqref{eq:alg2-square-norm-eta} gives
	\begin{dmath}
		\label{eq:alg2-sum-square-norm-b-eta}
		\rho_s \left ( \left | \bar{b}^{n-1} \right |^2-\left | \bar{b}^n \right |^2 \right ) + \rho_r \left ( \left | \bar{\eta}^{n-1} \right |^2 - \left | \bar{\eta}^n \right |^2 \right )= - \rho_s \rho_r \left ( \rho_r + \rho_s s^2 \right ) \left | \mu^n - \eta^{n-1} \right |^2 
		- \rho_s \rho_r \left (  \rho_s + \rho_r r^2 \right ) \left | b^{n-1} - q^n \right |^2 + 2 \rho_s \rho_r \left \langle \bar{b}^{n-1}, \eta^{n-1} - \mu^n + r \left ( b^{n-1} - q^n \right) \right \rangle 
		- 2 \rho_r \rho_s \left \langle \bar{\eta}^{n-1}, b^{n-1} - q^n -s \left ( \eta^{n-1} -  \mu^n \right) \right \rangle\\
		+ 2 \rho_r \rho_s \left ( \rho_r r - \rho_s s  \right ) \left \langle \mu^n - \eta^{n-1}, b^{n-1} - q^n \right \rangle.
	\end{dmath}
	Now, let us deal with the last terms in \eqref{eq:alg2-sum-square-norm-p-nu} and \eqref{eq:alg2-sum-square-norm-b-eta}. Using the Cauchy-Shwartz inequality, we get that 
	\begin{dmath}
		\label{eq:alg2-cauchy-shwartz-p-nu}
		2 \rho_r \rho_s \left ( \rho_r r - \rho_s s  \right ) \left \langle \mu^n - \nu^{n-1}, B\phi^n - p^{n-1} \right \rangle \geq
		- \rho_r \rho_s \left | \rho_r r - \rho_s s  \right | \left ( \alpha \left | \mu^n - \nu^{n-1} \right |^2 + \frac{1}{\alpha} \left | B\phi^n - p^{n-1} \right |^2 \right ), \forall \alpha > 0,
	\end{dmath}
	and
	\begin{dmath}
		\label{eq:alg2-cauchy-shwartz-b-eta}
		2 \rho_r \rho_s \left ( \rho_r r - \rho_s s  \right ) \left \langle \mu^n - \eta^{n-1}, b^{n-1} - q^n \right \rangle \geq 
		- \rho_r \rho_s \left | \rho_r r - \rho_s s  \right | \left ( \beta \left | \mu^n - \eta^{n-1} \right |^2 + \frac{1}{\beta}\left |  b^{n-1} - q^n \right |^2 \right ), \forall \beta > 0.
	\end{dmath}
	Considering \eqref{eq:alg2-cauchy-shwartz-p-nu} and \eqref{eq:alg2-sum-square-norm-p-nu}, we have
	\begin{dmath}
		\label{eq:alg2-sum-square-norm-p-nu-cauchy-shwartz}
		\rho_s \left ( \left | \bar{p}^{n-1} \right |^2-\left | \bar{p}^n \right |^2 \right ) + \rho_r \left ( \left | \bar{\nu}^{n-1} \right |^2-\left | \bar{\nu}^n \right |^2 \right ) \geq 
		- \rho_r \rho_s \left ( \rho_r + \rho_s s^2 + \alpha \left | \rho_r r - \rho_s s  \right | \right ) \left | \mu^n - \nu^{n-1} \right|^2 
		- \rho_r \rho_s \left ( \rho_r r^2 + \rho_s +\frac{1}{\alpha} \left | \rho_r r - \rho_s s  \right | \right ) \left | B \phi^n - p^{n-1}  \right |^2 
		+ 2 \rho_s \rho_r \left \langle \bar{p}^{n-1}, \mu^n - \nu^{n-1} + r \left ( p^{n-1} - B\phi^n \right) \right \rangle 
		- 2 \rho_r \rho_s \left \langle \bar{\nu}^{n-1}, B\phi^n - p^{n-1} -s \left ( \nu^{n-1} -  \mu^n \right) \right \rangle.
	\end{dmath}
	Similarly for \eqref{eq:alg2-cauchy-shwartz-b-eta} and \eqref{eq:alg2-sum-square-norm-b-eta}, we get
	\begin{dmath}
		\label{eq:alg2-sum-square-norm-b-eta-cauchy-shwartz}
		\rho_s \left ( \left | \bar{b}^{n-1} \right |^2-\left | \bar{b}^n \right |^2 \right ) + \rho_r \left ( \left | \bar{\eta}^{n-1} \right |^2 - \left | \bar{\eta}^n \right |^2 \right ) \geq
		- \rho_s \rho_r \left ( \rho_r + \rho_s s^2 + \beta \left | \rho_r r - \rho_s s  \right | \right ) \left | \mu^n - \eta^{n-1} \right |^2 
		- \rho_s \rho_r \left (  \rho_s + \rho_r r^2 + \frac{1}{\beta} \left | \rho_r r - \rho_s s  \right | \right ) \left | b^{n-1} - q^n \right |^2 
		+ 2 \rho_s \rho_r \left \langle \bar{b}^{n-1}, \eta^{n-1} - \mu^n + r \left ( b^{n-1} - q^n \right) \right \rangle 
		- 2 \rho_r \rho_s \left \langle \bar{\eta}^{n-1}, b^{n-1} - q^n -s \left ( \eta^{n-1} -  \mu^n \right) \right \rangle.
	\end{dmath}
	Summing inequalities \eqref{eq:alg2-sum-all-extremal-ineq}, \eqref{eq:alg2-sum-square-norm-p-nu-cauchy-shwartz}, \eqref{eq:alg2-sum-square-norm-b-eta-cauchy-shwartz}, and choosing $\alpha= \beta= 1$, we obtain
	\begin{dmath}
		\label{eq:alg2-final-ineq}
		\rho_s \left ( \left | \bar{p}^{n-1} \right |^2-\left | \bar{p}^n \right |^2 + \left | \bar{b}^{n-1} \right |^2-\left | \bar{b}^n \right |^2 \right ) 
		+ \rho_r \left ( \left | \bar{\nu}^{n-1} \right |^2-\left | \bar{\nu}^n \right |^2 + \left | \bar{\eta}^{n-1} \right |^2-\left | \bar{\eta}^n \right |^2\right ) \geq 
		\rho_r \rho_s \left( 2  s - \rho_r - \rho_s s^2 - \left | \rho_r r - \rho_s s  \right | \right) \left ( \left | \mu^n - \nu^{n-1} \right|^2 + \left | \mu^n - \eta^{n-1} \right|^2  \right ) 
		+ \rho_r \rho_s \left ( 2 r - \rho_r r^2 - \rho_s - \left | \rho_r r - \rho_s s  \right | \right ) \left (  \left | B \phi^n - p^{n-1}  \right |^2 + \left | b^{n-1} - q^n  \right |^2 \right ) 
		+2 \rho_r \rho_s \left( I\left ( \mu^n \right) - I \left( \mu^* \right) - \left \langle \bar{\mu}^n, p^* - b^* \right \rangle \right ).
	\end{dmath}
	It follows from \eqref{eq:alg2-final-ineq}, that for $\rho_r, \rho_s >0$ satisfying the following system of inequalities
	\begin{equation}
	\left\{\begin{array}{lr}
	2  s - \rho_r - \rho_s s^2 - \left | \rho_r r - \rho_s s  \right | & >0, \\
	2 r - \rho_r r^2 - \rho_s - \left | \rho_r r - \rho_s s  \right | & >0,
	\end{array}\right.
	\end{equation}
	the sequence 
	\[
	\rho_s \left ( \left | \bar{p}^{n} \right |^2  + \left | \bar{b}^n \right |^2 \right ) + \rho_r \left ( \left | \bar{\nu}^n \right |^2 + \left | \bar{\eta}^n \right |^2 \right )
	\]
	is decreasing, and the right hand in \eqref{eq:alg2-final-ineq} converges to zero (remember that $I\left ( \mu^n \right) - I \left( \mu^* \right)- \left \langle \bar{\mu}^n, p^* - b^* \right \rangle \geq 0 $).
	
	In view of Step 2, 3, 5 and 6 of Algorithm 4, we deduce that
	\begin{equation}
	\left | \bar{p}^{n} \right |^2,  \left | \bar{b}^n \right |^2, \left | \bar{\nu}^n \right |^2, \left | \bar{\eta}^n \right |^2, \left | B\bar{\phi}^{n} \right |^2,  \left | \bar{q}^n \right |^2, \left | \bar{\mu}^{n} \right |^2 \text{are bounded},
	\end{equation}
	and
	\begin{dmath}
		\lim_{n \to \infty} \left | B \phi^n - p^n \right | = \lim_{n \to \infty} \left | b^n - q^n \right | = \lim_{n \to \infty} \left |\mu^n - \nu^n \right | = \lim_{n \to \infty} \left |\mu^n - \eta^n \right | 
		= \lim_{n \to \infty} \left |\nu^n - \nu^{n-1} \right | = \lim_{n \to \infty} \left |\eta^n - \eta^{n-1} \right | =  \lim_{n \to \infty} \left | p^n - p^{n-1} \right | 
		=  \lim_{n \to \infty} \left | b^n - b^{n-1} \right | = \lim_{n \to \infty}  I\left ( \mu^n \right) - I \left( \mu^* \right) - \left \langle \bar{\mu}^n, p^* - b^* \right \rangle=0.
	\end{dmath}
	The statements in the theorem then follow from a similar argument to the one in the proof of the convergence of Algorithm 1.
\end{proof}

\ifCLASSOPTIONcaptionsoff
  \newpage
\fi



\bibliographystyle{IEEEtran}
%

%



